\documentclass{gtart}     
\usepackage{amssymb} 
\usepackage{epsfig}
\usepackage{float}
\usepackage{amsmath}
\usepackage{amsthm}
\usepackage{latexsym}
\usepackage{graphicx}  
\usepackage{renzo}
\usepackage{amsmath,amscd} 
\begin{document}
%
%
%
%
%
%
%
%
%

\title{Counting Bitangents with Stable Maps}                    
\authors{David Ayala and Renzo Cavalieri}                  
\address{Department of Mathematics\\ Stanford University\\
450 Serra Mall, Bldg. 380\\
Stanford, CA 94305-2125}
\secondaddress{Department of Mathematics\\
University of Utah\\
155 South 1400 East, Room 233\\
Salt Lake City, UT 84112-0090}                  


\email{ayala@math.Stanford.EDU\\renzo@math.utah.edu}                     
%
 
\begin{abstract}   

This paper is an elementary introduction to the theory of moduli spaces of curves and maps. As an application to enumerative geometry, we show how to count the number of bitangent lines to a projective plane curve of degree $d$ by doing intersection theory on moduli spaces.

\end{abstract}



%
%
%
\makeshorttitle  

%
\section*{Introduction}

\subsection*{Philosophy and Motivation}
The most apparent goal of this paper is to answer the following enumerative question:

\vspace{0.5cm}
\textit{``What is the number $N_{\mathcal{B}}(d)$ of bitangent lines to a generic projective plane curve $Z$ of degree d?''}
\vspace{0.5cm}

This is a very classical question, that has  been successfully solved with fairly elementary methods (see for example~\cite{gh:ag}, page 277). Here we propose to approach it from a very modern and 
``technological'' angle: we think of lines in the projective plane as
maps $\mu:\proj\rightarrow\Proj$ of degree $1$. We mark two points
$p_1$ and $p_2$ on $\proj$ and keep track of their image via the map
$\mu$. We then construct the space of all such marked maps, and ask
ourselves: can we understand the locus $\mathcal{B}$ of all maps that
are tangent to $Z$ at the images of both $p_1$ and $p_2$? The answer
fortunately is yes. The description of $\mathcal{B}$ allows us to produce 
a closed formula for $N_{\mathcal{B}}$ in all degrees.

This brief description already reveals that there is something deep
and interesting going on here, and that the journey is much more
important than the destination itself. Our major goal is to introduce
the reader to the rich and beautiful theory of Moduli Spaces in a
hopefully ``soft'' way, with the final treat of seeing it concretely
applied to solve our classical problem. 

It is our intention for this paper to be a very readable expository
work. We designed it to be accessible to a first year graduate student who is
considering algebraic geometry as a specialty field. We emphasize
geometric intuition and visualization 
above all, at the cost of silently glossing over some
technical details here and there.

\subsection*{Outline of the Paper}
This paper is divided into three sections, getting progressively more
advanced.
\vspace{0.5cm}

The first section introduces some basic ideas and techniques in modern
algebraic geometry, necessary to develop and understand the later two
sections and is intended for the unexperienced reader. We quickly
tread through the most basic ideas in intersection theory; we
introduce the concept of families of algebro-geometric objects; we
discuss the specific example of vector bundles, and give a working
sketch of the theory of Chern classes. Finally, we describe two
interesting constructions: the blow-up and jet bundles.\\ 
Entire books
have been written on each of these topics, so we have no hope or
pretense to be complete, or even accurate. Yet, we still think it
 valuable to present what lies in the back of a working mathematician's mind, in the firm belief that a solid geometric intuition is the best stairway to understand and motivate the technicalities and abstract generalizations needed to make algebraic geometry ``honest''.

\vspace{0.5cm}
The second section is the development of most of the theory. After a quick qualitative introduction to  moduli spaces, we  discuss our  main characters: the moduli spaces of rational  stable curves, and of rational stable maps. Intersection theory
on the moduli spaces of  stable maps, commonly referred to as Gromov-Witten Theory, is currently an extremely active area of
research.

\vspace{0.5cm} Finally, in the third section we apply all the theory
developed so far to solve the bitangent problem. We explore in further
detail the moduli spaces of rational stable maps of degree $1$ to
$\Proj$, with one and two marked points. By intersecting appropriate cycle classes on these spaces we extract one of the classical \textit{Pl\"ucker Formulas}, expressing the number of bitangents as a function of the degree $d$ of the curve.

\subsection*{References}

We suggest here some canonical references for the reader in search of
more rigor and completeness. For intersection theory,~\cite{f:it} is
a fairly technical book, but definitely it has the last word on it. It
also presents Chern classes from an algebraic point of view. A
discussion of Chern classes from a geometric point of view can be
found in~\cite{bt:df}.\\ A good treatment of blow-ups can be found in
any basic book in algebraic geometry, for example~\cite{gh:ag} or
~\cite{h:ag}.\\ A very pleasant reference for jet bundles is
~\cite{v:jb}. An extensive treatment of jet bundles is found in ~\cite{s:jb}.\\ Our presentation of moduli spaces follows the spirit of
~\cite{k:kf}; for anybody interested in getting serious,~\cite{hm:mc} is the way to go. Finally, a good introduction to $\psi$ classes is~\cite{k:pc}.

\subsection*{Acknowledgments}
We first of all owe the inspiration for this work to Joachim Kock, who
outlines this strategy for counting bitangents in his talk~\cite{k:bt}, and is also
responsible, with Israel Vainsencher, for the best elementary introduction to Gromov-Witten
Theory we know of,~\cite{k:kf} . We also  thank Aaron Bertram, Herb Clemens, Tommaso de Fernex, Hugo Rossi and Ravi Vakil for their comments, suggestions and encouragement.

\section{Preliminaries}
\subsection{Intersection Theory}
It will be helpful, but not essential, that the reader be familiar with the Chow ring,
$A^\ast(X)$, of an algebraic variety $X$.  The ring\footnote{This is probably our greatest sloppiness. In order for  $A^\ast(X)$ to be a ring we need $X$ to be smooth. Since the spaces we will actually work with satisfy these hypotheses, we do not feel too guilty.} $A^\ast(X)$ is, in some loose sense,  the algebraic
counterpart of the cohomology ring $H^\ast(X)$, and it allows us to make precise in the algebraic
 category the intuitive concepts of oriented intersection 
in topology.

We think of elements of the group $A^n(X)$ as formal finite sums of
codimension $n$ closed subvarieties (cycles), modulo an equivalence relation called rational equivalence.
  $A^\ast(X)=\oplus_0^{dim X}A^n(X)$ is a graded  ring with product given by intersection.
\\ Intersection
is
independent of the choice of representatives for the equivalence classes.

In topology, if
we are interested in the cup product of  two cohomology classes $\mathbf{a}$ and $\mathbf{b}$,
we can choose  representatives $a$ and $b$ 
that are transverse to each other. We can assume this since transversality is a
generic condition: if $a$ and $b$ are not transverse then
we can perturb them ever so slightly and make them transverse while
not changing their classes. This
being the case, then $a \cap b$ represents the cup product class $\mathbf{a} \cup \mathbf{b}$.

In algebraic geometry, even though this idea must remain the backbone of our intuition, things
are a bit trickier. We will soon see examples of cycles that are rigid, in the sense
that their representative is unique, and hence ``unwigglable". Transversality then becomes an unattainable
dream. Still, with the help of substantially sophisticated machinery (the interested reader can
consult~\cite{f:it}), we can define an algebraic version of intersection classes and  a product
 that reduces
to the ``geometric" one when transversality can be achieved.
 
Throughout this paper, a bolded symbol will represent a class, the unbolded symbol a geometric representative. The intersection of two classes $\mathbf{a}$ and \textbf{b} will be denoted by \textbf{ab}.

\vspace{0.5cm}

\textbf{Example}: the Chow Ring of Projective Space. 
$$A^\ast(\mathbb{P}^n) = \frac{\mathbb{C}[\mathbf{H}]}{(\mathbf{H}^{n+1})},$$ where
$\mathbf{H}\in A^1(\mathbb{P}^n)$ is the class of a hyperplane
$H$. 

\subsection{Families and Bundles}

One of the major leaps in modern algebraic geometry comes from the insight that, to fully understand  algebraic varieties, we should
not study them one by one, but understand how they  organize themselves in 
families.\\
We are all familiar, maybe subconsciously, with the concept of a family. When, in high
school, we dealt with ``all parabolas of the form $y=ax^2$" or ``all circles with
center at the origin", we had in hand prime examples of families of algebraic varieties.\\
The idea is quite simple: we have a parameter space, $B$, called the base of the
family. For each point $b\in B$ we want an algebraic variety $X_b$ with certain properties.
Further, we want all such varieties to be organized together to form an algebraic variety $\mathcal{E}$,
 called the
total space of the family.

A little more formally we could define a family of objects of type $\mathcal{P}$
as a morphism of algebraic varieties
 $$
\begin{array}{c}
\mathcal{E}\\
\pi\downarrow\\
B
\end{array}
$$
where $\pi^{-1}(b)$  is an object of type $\mathcal{P}$. 


A \textbf{section} of a family $\pi:\mathcal{E}\rightarrow B$ is a map $s:B\rightarrow \mathcal{E}$
such that $\pi\circ s:B\rightarrow B$ is the identity map. Often,
the section $s$ is written
$$
\begin{array}{c}
\mathcal{E}\\
\pi\downarrow\uparrow s\\
B.
\end{array}
$$
Notice that $s(b)\in\pi^{-1}(b)$. 

\vspace{0.5cm}
Given a family
$\pi:
\mathcal{E}\longrightarrow B$ and a map $f:M\rightarrow B$ we can construct a new family
$$
\begin{array}{c}
f^\ast \mathcal{E}\\
\downarrow\pi_f\\
M,
\end{array}
$$
called the \textbf{pull-back of $\pi$ via $f$}:
$$f^{\ast}\mathcal{E}=\{(m,e)\in M\times \mathcal{E} \mid f(m)=\pi(e)\}.$$
 Intuitively, the fibre of $\pi_f$ over 
a point $m\in M$ will be the fibre of $\pi$ over $f(m)$. An essential property of
this construction is that it is natural, up to isomorphism.

\subsubsection{Vector Bundles}
 
A \textbf{vector bundle} of rank $n$ is a family
$\pi:\mathcal{E}\rightarrow B$ of vector spaces over $\mathbb{C}$ of dimension
$n$ which is \textbf{locally
trivial}\footnote{To be precise, more structure is needed: the clutching functions must take values in $GL(n,\mathbb{C})$.}. By locally trivial we mean that there is on open cover
$\{U_\alpha\}$ of $B$ such that
$\pi^{-1}(U_\alpha)\cong U_\alpha\times\mathbb{C}^n$. Our vector bundle is uniquely determined by how these trivial pieces
glue together. 

A vector bundle of rank one is called a \textbf{line bundle}, as its fibers are (complex) lines.

\vspace{0.5cm}

Given two vector bundles
$$
\begin{array}{ccc}
\mathcal{E}_1 & & \mathcal{E}_2\\
\pi_1\downarrow & $and$ & \pi_2\downarrow\\
B & & B
\end{array}
$$
over the same base space, one can define their Whitney sum
$$
\begin{array}{c}
\mathcal{E}_1\oplus \mathcal{E}_2\\
\pi\downarrow\\
B,
\end{array}
$$
where a fibre $\pi^{-1}(b)$ is the direct sum of the vector spaces
$\pi_1^{-1}(b)\oplus\pi_2^{-1}(b)$. It can be easily verified that
this  family satisfies the local triviality
condition.

Similarly, one can define the
tensor product $\mathcal{E}_1\otimes \mathcal{E}_2$, the dual bundle $\mathcal{E}^\ast$, the wedge product $\bigwedge^p(\mathcal{E})$  
and the
bundle $Hom(\mathcal{E}_1, \mathcal{E}_2)=(\mathcal{E}_2\otimes \mathcal{E}_1^\ast)$. 

\subsubsection{Characteristic Classes of Bundles}\label{Chern}
For every vector bundle  
there is a natural section $s_0:B\rightarrow
\mathcal{E}$ defined by $$s_0(b)=(b,0)\in\{b\}\times\mathbb{C}^n.$$  It is called
the zero section, and it gives an embedding of $B$ into $\mathcal{E}$.

A
natural question to ask is  if there exists another section $s:B\rightarrow \mathcal{E}$
which is disjoint form the zero section, i.e. $s(b)\not=s_0(b)$ for all
$b\in B$. 
The \textbf{Euler class} of this vector bundle ($\mathbf{e}(\mathcal{E})\in A^n(B)$)  is defined to be the class of the self-intersection of the zero section: it measures obstructions for the above question to be answered affirmatively. This means that $\mathbf{e}(\mathcal{E})=0$ if and only if a never vanishing section exists.
It easily follows from the Poincar\'e-Hopf theorem that for  a
manifold $M$, the following formula holds:
$$\mathbf{e}(TM)\cap[M]=\chi(M).$$ That is, the degree of the
Euler class of the tangent bundle is the Euler characteristic.

\vspace{0.5cm}
The Euler class of a vector bundle is the first and most important example of a whole family
of ``special" cohomology classes associated to a bundle, called the \textbf{Chern classes} 
of $\mathcal{E}$. The $k$-th Chern class of $\mathcal{E}$, denoted $\mathbf{c}_k(\mathcal{E})$, lives in $A^k(B)$.
In the literature you can find a wealth of definitions for Chern classes, some more geometric,
dealing with obstructions to finding a certain number of linearly independent sections of the bundle,
some purely algebraic. Such formal definitions, as important as they are (because they assure us that we
are talking about something that actually exists!), are not particularly illuminating. In concrete terms, what you really need to know is that Chern classes
are cohomology classes associated to a vector bundle that satisfy a series of really nice properties, which 
we are about to recall. 

Let $\mathcal{E}$ be a vector bundle of rank $n$:
\begin{description}
	\item [identity:] by definition, $\mathbf{c}_0(\mathcal{E}) = 1$.
	\item [normalization:] the $n$-th Chern class of $\mathcal{E}$ is the Euler class: $$\mathbf{c}_n(\mathcal{E})=\mathbf{e}(\mathcal{E}).$$
	\item [vanishing:] for all $k>n$,  $\mathbf{c}_k(\mathcal{E})=0$.
	\item [pull-back:] Chern classes commute with pull-backs:
$$f^\ast\mathbf{c}_k(\mathcal{E})=\mathbf{c}_k(f^\ast \mathcal{E}).$$
	  \item [tensor products:] if $L_1$ and $L_2$ are line bundles,
$$c_1(L_1\otimes L_2)= c_1(L_1) + c_1(L_2).$$
	\item [Whitney formula:] for every extension of bundles\label{whit}
$$ 0\rightarrow \mathcal{E}' \rightarrow \mathcal{E} \rightarrow \mathcal{E}'' \rightarrow 0 ,$$
 the $k$-th Chern class of $\mathcal{E}$ can be computed in terms of the Chern classes of $\mathcal{E}'$ and $\mathcal{E}''$,
by the following formula:
$$\mathbf{c}_k(\mathcal{E})=\sum_{i+j=k} \mathbf{c}_i(\mathcal{E}')\mathbf{c}_j(\mathcal{E}'').$$
\end{description}

Using the above properties it is immediate to see:
\begin{enumerate} 
	\item all the Chern classes of a trivial bundle vanish (except the $0$-th, of course);
	\item for a line bundle $L$, $\mathbf{c}_1(L^\ast)= -\mathbf{c}_1(L)$.
\end{enumerate}
\vspace{0.5cm}
To show how to use these properties   to work with Chern classes, we will
 now calculate the  first Chern class of the tautological line bundle over $\proj$. 
The tautological line bundle is  $$
\begin{array}{c}
\mathcal{S}\\
 \pi\downarrow \ \ \\
 \ \ \proj,
\end{array}
$$
where $\mathcal{S}=\{(p,l)\in\mathbb{C}^2\times\proj\mid p\in l\}$. It is called tautological because the fiber over a point in $\proj$ is the line  that point represents.

Our tautological family fits into the short exact
sequence of vector bundles over $\proj$
$$
\begin{array}{ccccccccc}
0 &\rightarrow & \mathcal{S} & \rightarrow & \mathbb{C}^2\times\proj & \rightarrow
& \mathcal{Q} & \rightarrow & 0 \\
  &            &            & \searrow     & \downarrow            & \swarrow & 
 & & \\
 & & & & \proj & & & & 
\end{array}         
$$
where $\mathcal{Q}$ is the 
bundle whose
fibre 
over a line $l\in\proj$ is the quotient vector space
$\mathbb{C}^2/l$. Notice that $\mathcal{Q}$
is also a line bundle. From the above sequence, we have that
\begin{eqnarray}
0=\mathbf{c}_1(\mathbb{C}^2\times\proj)=\mathbf{c}_1(\mathcal{S})+\mathbf{c}_1(\mathcal{Q}).\label{triv}
\end{eqnarray}

Since $\proj$ is topologically a  sphere,  which has Euler
characteristic 2, then
\begin{eqnarray}
2=\mathbf{c}_1(T\proj)=\mathbf{c}_1(\mathcal{S}^\ast)+\mathbf{c}_1(\mathcal{Q})
=-\mathbf{c}_1(\mathcal{S})+\mathbf{c}_1(\mathcal{Q}).\label{tang}
\end{eqnarray}
The second equality in \ref{tang} holds because
$T\proj$ is the line bundle $Hom(\mathcal{S},\mathcal{Q})= 
\mathcal{S}^\ast\otimes\mathcal{Q}$. 

It now follows from (\ref{triv}) and (\ref{tang}) that $\mathbf{c}_1(\mathcal{S})=-1$.

\subsection{Blow-up}

Let us begin by discussing the prototypical example of the blow up of a point on a surface:
first off, the blow up is a local
construction and so we need only understand the local picture. 
\\ Consider the map 
$$
\begin{array}{rccc}
\phi: & \mathbb{C}^2 & \rightarrow & \mathbb{C}\\
      & (x,y)        & \mapsto     & y/x.\\		
\end{array}
$$
This is a rational map and is not defined on
the line $\{x=0\}$. We may try to fix this by modifying our target space to
$\proj$. Still, $\phi$ cannot be defined at $\mathbf{0}:=(0,0)$. In fact,
 along any line $l$ through the origin, the limit of $\phi$ at $\mathbf{0}$ 
is the slope of $l$.    We
would like to modify $\mathbb{C}^2$ to a smooth surface 
birational to it, where the map $\phi$ can be defined everywhere. 
We would like points outside $\mathbf{0}$ to remain ``untouched" and
$\mathbf{0}$ to be replaced by a $\proj$ whose points represent all
tangent directions at $\mathbf{0}$. 

Here is how to do it: consider the graph of $\phi$,
$\Gamma_\phi\subset\mathbb{C}^2\times\proj$. We have the commutative
diagram
$$
\begin{array}{rrc}
 &  \Gamma_\phi & \subset\mathbb{C}^2\times\proj \\
 & \stackrel{(id,\phi)}{}{\nearrow}\swarrow \pi_1  &\downarrow \pi_2  \\
&\mathbb{C}^2\setminus\{0\}  \stackrel{\phi}{\longrightarrow} & \proj 
\end{array}
$$
 The closure
$\overline{\Gamma}_\phi$ is
 what we are looking for. It is birational to
$\mathbb{C}^2$; the left projection ${\pi_1}_{\mid\Gamma_\phi}$ is an isomorphism onto
$\mathbb{C}^2-\{\mathbf{0}\}$ and $\pi_1^{-1}(\mathbf{0})=\proj$. We define 
 the blow-up of $\mathbb{C}^2$ at  $\mathbf{0}$ as  $Bl_\mathbf{0}(\mathbb{C}^2):=\overline{\Gamma}_\phi=\Gamma_\phi\cup\proj$;
the projective line $\pi_1^{-1}(\mathbf{0})$ is called the exceptional divisor of the blow-up and denoted $E$. We
have  obtained a new (smooth!) space by replacing, in a particularly favorable way, a point
with the projectivization of its normal
bundle. 
\begin{figure}[htbp]
\begin{center}
\includegraphics[height=6cm]{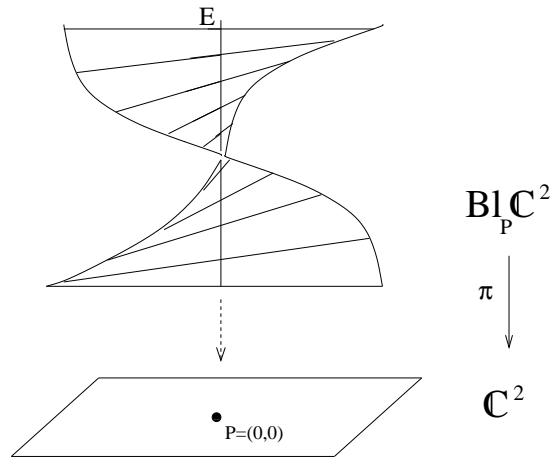}
\caption{the blow-up of $\mathbb{C}^2$ at the origin.}
\end{center}
\end{figure}

In general, for $Y\subset X$ a closed subvariety of codimension $k\geq0$,
one can construct a new space $Bl_Y(X)$ such that:
\begin{enumerate}
	\item $Bl_Y(X)$ is birational to $X$;
	\item points outside $Y$ are
untouched;
	\item  a point in $Y$ is replaced by
$\mathbb{P}^{k-1}$, representing the ``normal'' directions to
$Y$ at that point. 
\end{enumerate}

The total space of the blow-up of $\mathbb{C}^2$ at $\mathbf{0}$ 
admits a natural map to the exceptional divisor, consisting of projecting points along
lines through the origin. This realizes $Bl_0 (\mathbb{C}^2)\rightarrow E$ as a line
bundle over $\proj$. This is the tautological bundle, which does not have any global sections besides the $0$-section. It follows  that the class $\mathbf{E}$
of the exceptional divisor admits only one representative, namely $E$ itself.
It is therefore impossible to compute the self-intersection $\mathbf{EE}$ by means of
intersecting two transverse representatives of the class.

\subsection{Jet Bundles}\label{jets}
Let
$L$ be a line bundle over a variety, $X$. Then the local sections of this line bundle
form a vector space. In fact, locally, such a
section is just a complex valued function on some open set in X. 
We will now describe a new
vector bundle over $X$ whose fibre over $x\in X$ consists of all Taylor
expansions of these sections centered about $x$ and truncated after
degree $k$. To
see how the locally
trivial charts of this bundle glue together is simply a matter
of shifting the center
of a Taylor expansion. We call this bundle the $k$th jet bundle of $L$
and denote it by $J^kL$. 
In particular $J^0L=L$.\\
Notice that the first jet bundle keeps track
of all locally defined functions and differential forms and so there is an
obvious surjection $J^1L\rightarrow L$. This  gives us the short
exact sequence
$$
0\rightarrow L\otimes\Omega\rightarrow J^1L\rightarrow
L\rightarrow0
$$
which will be an essential tool later on.

The previous statement is a particular case of what can be considered the ``fundamental
theorem of jet bundles".

\newtheorem{teor}{Theorem}

\begin{teor}
For all $n\geq 0$, the sequence
\begin{eqnarray}
0\rightarrow L\otimes Sym^n\Omega\rightarrow J^nL\rightarrow
J^{n-1}L\rightarrow0
\label{jet}
\end{eqnarray}
is exact.
\end{teor}
For a slightly more rigorous  and still enjoyable account of jet bundles, refer to~\cite{v:jb}.

\section{Moduli Spaces}
\subsection{ A ``High School" Example}

What is the idea of a moduli space? A moduli space  of geometric objects of
a certain type is a
space which ``encodes'' information about collections of geometric objects of a given type, in the sense that:
\begin{enumerate}
	\item  points in the moduli space correspond bijectively to the desired geometric objects;
	\item  the moduli space itself has an algebraic structure that respects how the objects can organize
themselves in families.
\end{enumerate}

 For
example, suppose
that we would like to consider the space of all circles in the
plane.  Since a circle is uniquely the zero locus of a second
degree polynomial
of the form $(x-x_0)^2+(y-y_0)^2-r^2$,  upon
specifying the coordinates of the center  and its radius, we
have completely identified the circle. Thus, the space of all circles
in the plane can be represented by $\mathcal{M}:=\mathbb{R}^2\times\mathbb{R}_+$. This is indeed much more than just a set-theoretic correspondence. 

Consider
the tautological family

$$
\begin{array}{c}
\mathcal{U} \ \\
\downarrow\pi\\
\mathcal{M}
\end{array}
$$
where
$\mathcal{U}:=\{((x_0,y_0),r,(x,y))\mid (x-x_0)^2+(y-y_0)^2=r^2\}\subset\mathcal{M}\times\mathbb{R}^2$
and $\pi$ is the projection onto the first factor. This family
enjoys the following properties, that clarify the vague point $2$ above:
\begin{enumerate}
		\item for any family of
circles in the plane $p:E\rightarrow B$, there is a map
$m:B\rightarrow\mathcal{M}$ defined by $m(b)=p^{-1}(b)$;
	\item to every map $m:B\rightarrow\mathcal{M}$
 there uniquely corresponds a family of circles parametrized by $B$,
i.e. 
$$
\begin{array}{c}
m^\ast \mathcal{U}\\
p\downarrow\\
B
\end{array}$$
such that the fibre $p^{-1}(b)$ is the circle $m(b)$. 
\end{enumerate}
This is the best that we could have hoped for. In this case we say that $\mathcal{M}$
is a \textbf{fine moduli space} with $\mathcal{U}$ as its \textbf{universal family}.

Often, due to the presence of automorphisms of the parametrized objects, it is  impossible to achieve this perfect bijection between families of objects and morphisms to the moduli space. If only property 1 holds  we call the moduli space \textbf{coarse}.  

\subsection{Moduli of $n$ Points on $\proj$}

Let us now consider the moduli space $M_{0,n}$ of
all isomorphism classes of $n$ ordered distinct marked points
$p_i\in\proj$. The subscript 0 is to denote the genus of our curve
$\proj$. Since the  automorphism group $Aut(\proj)=PSL_2(\mathbb{C})$ allows 
us to move any three points on
$\proj$ to the ordered triple $(0,1,\infty)$, the space
$M_{0,n}$ reduces to a single point for $n\leq3$.

Going one step up,
$M_{0,4}=\proj-\{0,1,\infty\}$ : given a quadruple $(p_1,p_2,p_3,p_4)$, 
we can always perform the unique 
automorphism of $\proj$ sending $(p_1,p_2,p_3)$ to $(0,1,\infty)$; the isomorphism
class of the quadruple is then  determined by the image of the fourth point.

The general case is similar. Any $n$-tuple $\underline{p}=(p_1,\ldots, p_n)$ is equivalent to a $n$-tuple of the form $(0,1,\infty,\phi(p_4),\ldots,\phi(p_n))$, where $\phi$  is the unique automorphism of $\proj$ sending $(p_1,p_2,p_3)$ to $(0,1,\infty)$. This shows
$$M_{0,n}=\overbrace{M_{0,4}\times...\times M_{0,4}}^{n-3\  \mbox{times}} \setminus \{\mbox{all
diagonals}\}.$$ 
If we define
  $U_n:=M_{0,n}\times\proj$, then the projection of $U_n$ onto the first factor
  gives rise to a universal family
$$
\begin{array}{c}
U_n\\
\pi\downarrow\uparrow\sigma_i\\
M_{0,n}
\end{array}
$$
where the $\sigma_i$'s are the
 universal sections:
\begin{itemize}
	\item $\sigma_i(\underline{p})= (\underline{p},\phi(p_i))\in U_n$.
\end{itemize}
 This family is tautological since the fibre
 over a moduli point, which is the class of a marked curve, is the marked curve
 itself. \\
With $U_n$ as its universal family, 
$M_{0,n}$
 becomes a fine moduli space
  for isomorphism classes of  $n$ ordered distinct marked
  points on $\proj$.

This is all fine except $M_{0,n}$ is
not compact for $n\geq 4$. There are many reasons why compactness
is an extremely desirable property for moduli spaces. As an extremely practical
reason, proper (and if possible projective) varieties are much better behaved and understood
than non compact ones. Also, a compact moduli space encodes information on
how our objects can degenerate in families.
For example, what  happens
when $p_1\rightarrow p_2$
in $M_{0,4}$? 

In general there are many ways to compactify a space.
A ``good'' compactification $\overline{\mathcal{M}}$ of a moduli space $\mathcal{M}$ should have the following properties:
\begin{enumerate}
\item  $\overline{\mathcal{M}}$ should be itself a moduli space, parametrizing some natural generalization of the objects of $\mathcal{M}$.
\item  $\overline{\mathcal{M}}$ should not be a horribly singular space.
\item  the boundary $\overline{\mathcal{M}}\setminus\mathcal{M}$ should be a normal crossing divisor.
\item  it should be possible to describe boundary strata combinatorially in terms of simpler objects. This point may appear mysterious, but it will be clarified by the examples of stable curves and stable maps.
\end{enumerate}
 
In the case of rational $n$-pointed curves there is a definite winner among compactifications.

\subsection{Moduli of Rational Stable Curves}

We will discuss the simple example of $M_{0,4}$; this hopefully will, 
without submerging  us in combinatorial technicalities, provide intuition on the ideas and 
techniques used to
compactify the moduli spaces of $n$-pointed rational curves.

A  natural
first attempt would be to just allow the points to come together, i.e. enlarge
the collection of objects that we are considering from 
$\proj$ with $n$
ordered distinct marked points to $\proj$ with $n$ ordered,
not necessarily distinct, marked points. 

However, this will not quite work. For instance, consider the
families $$C_t=(0,1,\infty,t) \mbox{\ \  and\ \ }
D_t=(0,t^{-1},\infty,1).$$ For each $t\not=0$, up to an automorphism of
$\proj$, $C_t=D_t$, thus corresponding to the same
point in $M_{0,4}$. But for $t=0$, $C_0$ has $p_1=p_4$ whereas $D_0$
has $p_2=p_3$. These configurations are certainly not equivalent up to
an automorphism of $\proj$ and so should be considered as distinct
points in our compactification of $M_{0,4}$. Thus, we have a family with
two distinct limit points (in technical terms we say that the space is \textbf{nonseparated}).
This is not good. 

Our failed attempt was not completely worthless though since it allowed us to understand
that we want the condition $p_1=p_4$ to coincide with $p_2=p_3$, and likewise for
the other two possible disjoint pairs. On the one hand this is very promising: $3$ is the number of points needed to 
compactify $\proj\setminus\{0,1,\infty\}$ to $\proj$.
On the other hand, it is now  mysterious what modular interpretation to give to this compactification.

\begin{figure}[htbp]
\begin{center}
\includegraphics[height=6cm]{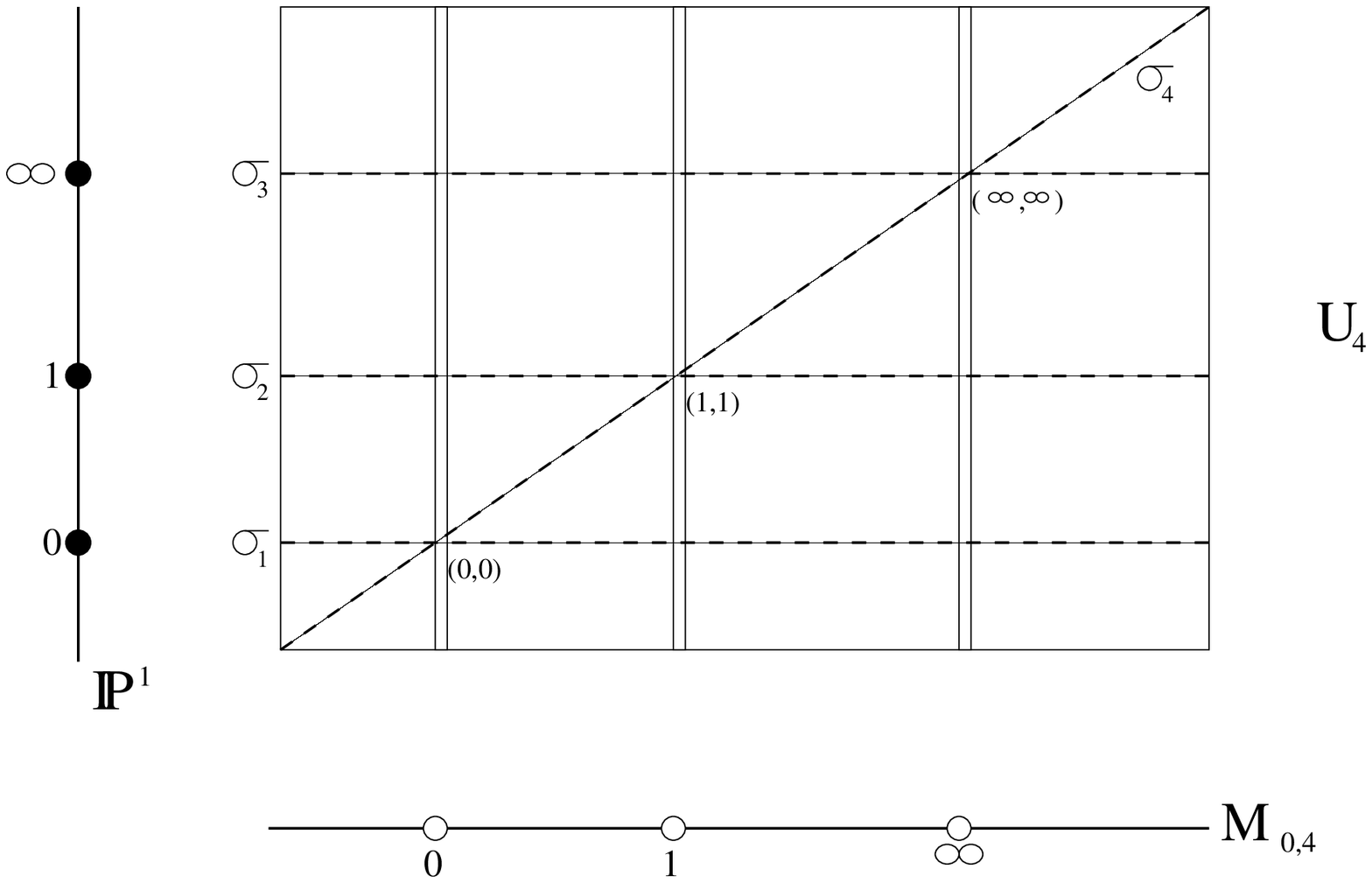}
\caption{first attempt at compactifying $U_{4}$}
\label{uf}
\end{center}
\end{figure}

To do so, let us turn carefully to our universal family, illustrated in Figure~\ref{uf}. The natural first step is to fill in the three points on the base, to complete
$U_4$ to $\proj\times\proj$ and extend the sections by continuity.

We immediately notice a bothersome asymmetry in this picture: the point $p_4$ is the only one allowed to come together with all the other points: yet common sense, backed up by  the explicit example just presented, suggests  that there should be democracy among the four points.
This fails where the diagonal section $\sigma_4$ intersects the three constant ones, i.e. at the
three points $(0,0)$, $(1,1)$, $(\infty,\infty)$. Let us blow-up $\proj\times\proj$ at these three points. This will
make all the sections disjoint, and still preserve the smoothness and projectivity of our universal
family. 
\\
The fibres over the three exceptional points are $\proj \cup E_i$:  nodal rational curves.
These are the new objects that we have to allow in order to obtain a good
compactification of $M_{0,4}$. 
\\
Let us finally put everything together, and state things carefully.

\newtheorem{defi}{Definition}

\begin{figure}[htbp]
\begin{center}
\includegraphics[height=4cm]{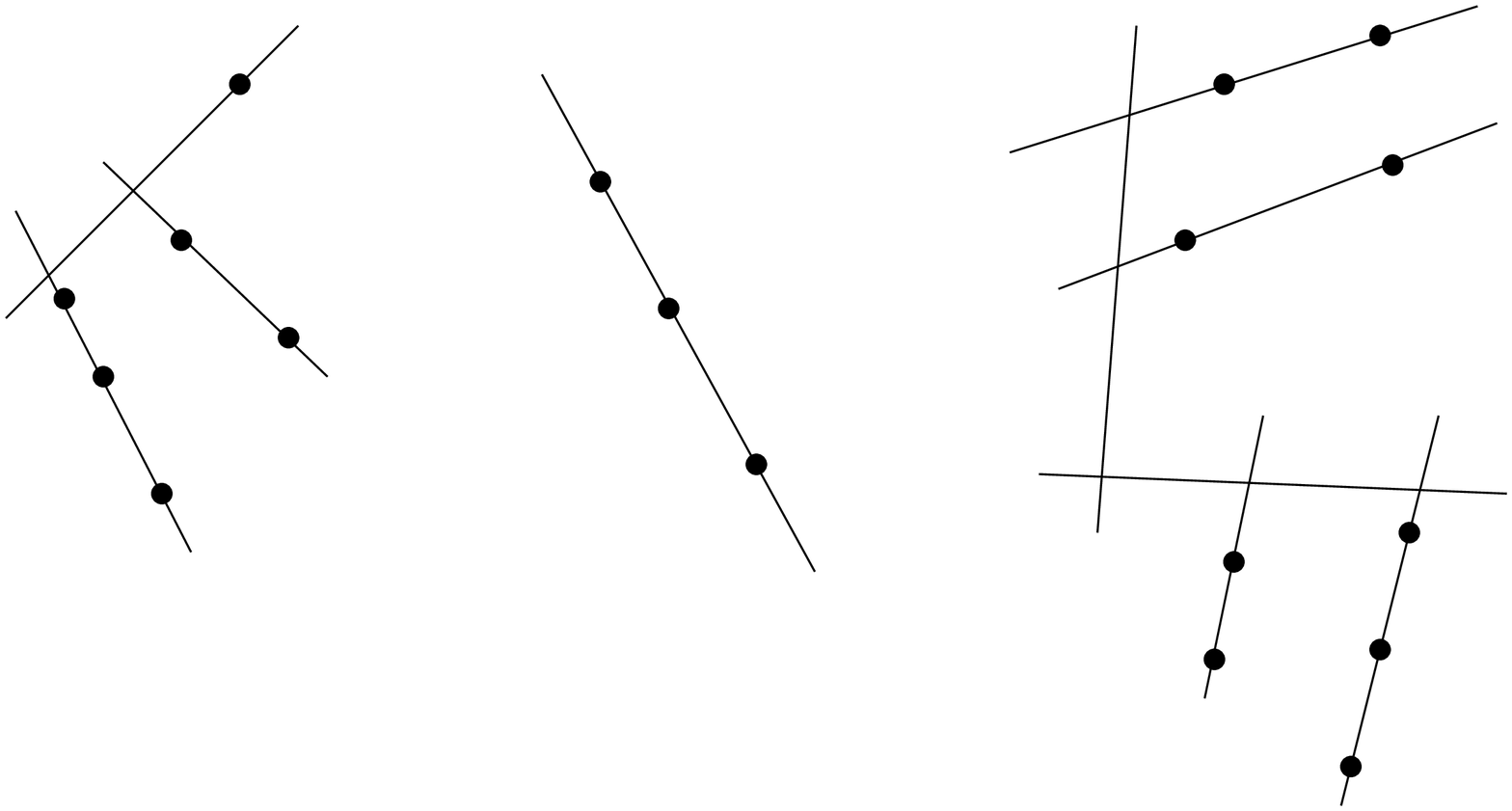}
\caption{stable marked trees.}
\label{mt}
\end{center}
\end{figure}

\begin{defi}
A \textbf{tree of projective lines} is a connected curve with the
  following properties: 
\begin{enumerate} 
	\item Each irreducible component is isomorphic
  to $\proj$. 
	\item The points of intersection of the
  components are ordinary double points. 
	\item There are no closed
  circuits, i.e., if a node is removed then the curve becomes
  disconnected.
\end{enumerate}
\end{defi}
 These three properties are equivalent to saying that
  the curve has arithmetic genus zero.  Each irreducible component
  will be called a twig. We will often draw a marked tree as
  in fig \ref{mt}, where each line represents a twig.

\begin{defi}
A marked tree is
\textbf{stable} if every twig has at least three special points (marks or
nodes). 
\end{defi}

This stability condition is
 equivalent to the existence of no nontrivial
automorphisms of the tree that fix all of the marks.

\begin{defi}
 
 $\overline{M}_{0,4}\cong \proj$ is the moduli space of isomorphism classes of four
pointed stable trees. It is a fine moduli space,
with universal family \linebreak $U_4 = Bl(\proj\times \proj)$.
\end{defi}
\vspace{0.5cm}

These
results generalize to larger $n$.\\
\textbf{Fact:}  The space $\overline{M}_{0,n}$ of
$n$-pointed rational stable curves is a fine moduli space compactifying
$M_{0,n}$. It is projective, and the universal family $\overline{U}_n$
is obtained from $U_n$ via a finite sequence of blow-ups. (In particular all the diagonals need to be blown up in an appropriate order) .  For further details see ~\cite{k:kf} or \cite{k:2}, \cite{k:1}.

\vspace{0.5cm}

 One of the exciting features of this theory is that all
these spaces are related to one another by natural morphisms.
Consider the map
$$\pi_i:\overline{M}_{0,n+1}\rightarrow\overline{M}_{0,n},$$
 defined by
forgetting the $i$th mark. It is obviously defined if the $i$th mark
does not belong to a twig with only three special points. If it does
belong to such a twig, then our resulting tree will no longer be
stable. In this case, we must perform what is called contraction. 
\begin{description}
	\item [Contraction:] We need to consider two cases:
	\begin{enumerate} 
		\item The remaining two special points are
both nodes. We make the tree again stable by contracting this twig so that the two nodes
are now one (see Figure \ref{c1}).

\begin{figure}[htbp]
\begin{center}
\includegraphics[width=13.5cm]{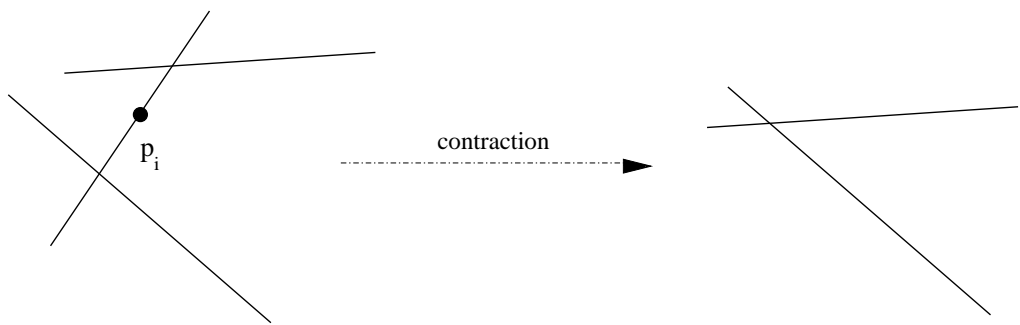}
\caption{contracting a twig with only two nodes.}
\label{c1}
\end{center}
\end{figure} 
		\item There is one
other mark and
one node on the twig in question. We make the tree stable by forgetting the twig and placing the mark
where the node used to be (Figure\ref{c2}). 

\begin{figure}[htbp]
\begin{center}
\includegraphics[width=14cm]{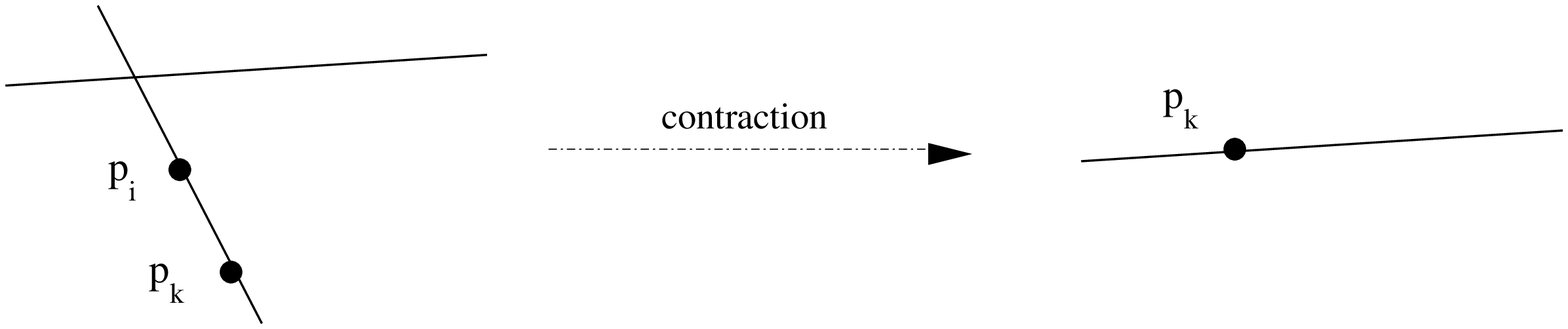}
\caption{contracting a twig with one node and one mark.}
\label{c2}
\end{center}
\end{figure}
	\end{enumerate}
\end{description}

\vspace{0.5cm}
We would like to construct a section 
$\sigma_i$ of the
family
$$
\begin{array}{c}
\overline{M}_{0,n+1}\\
\pi_k\downarrow\uparrow\sigma_i\\
\overline{M}_{0,n}
\end{array}
$$
 by defining the $k$th mark  to coincide with  $i$th one. It should trouble you 
that in doing so we are not considering curves with distinct
marked points, but we can get around this problem  by ``sprouting'' a new twig so that the node
is now where the $i$th mark was. The $k$th and the $i$th points now belong to this new twig.\\
This
process of
making stable a tree with two coinciding points is called \textbf{stabilization}.
\begin{figure}[htbp]
\begin{center}
\includegraphics[width=13cm]{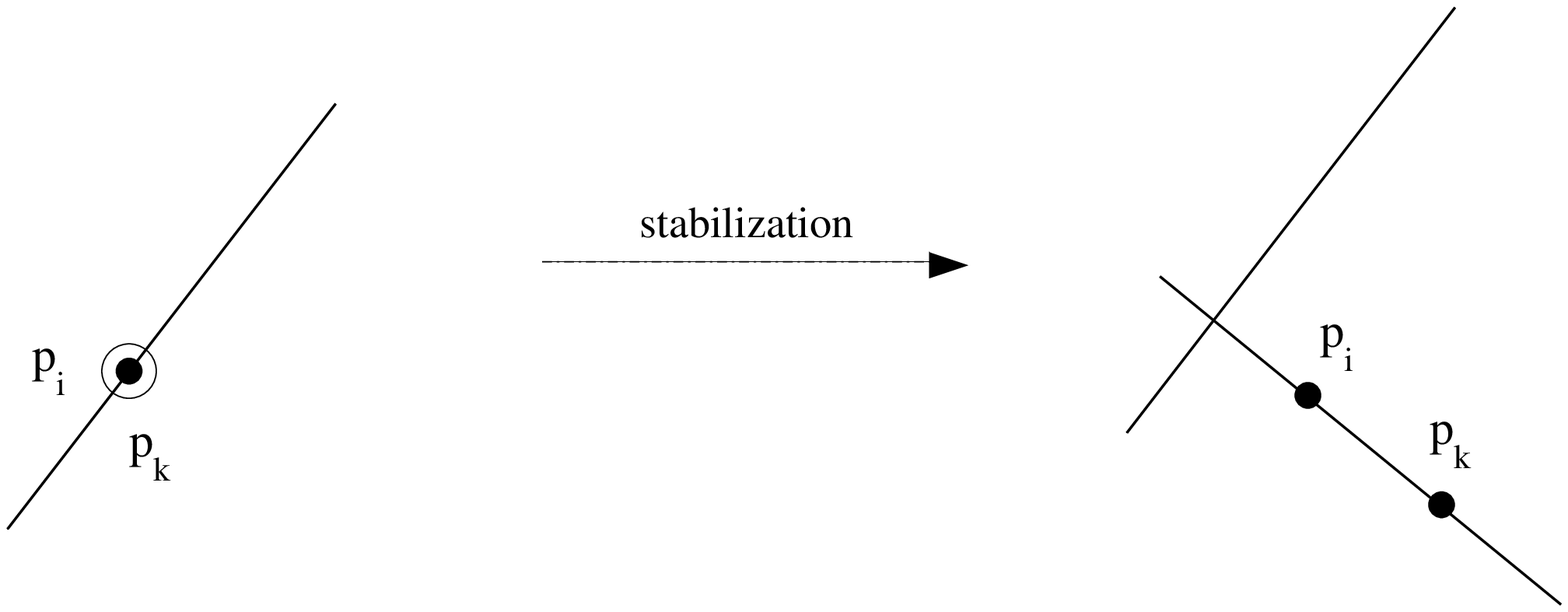}
\caption{stabilization}
\end{center}
\end{figure}

Finally, we may now identify our universal family
$$
\begin{array}{ccc}
\overline{U}_n & & \overline{M}_{0,n+1} \\
\pi\downarrow & $with the family$ & \pi_i\downarrow\\
\overline{M}_{0,n} & & \overline{M}_{0,n}
\end{array}
$$
as follows.\\
 The fibre
$\pi^{-1}([(C,p_1, \cdots, p_n)])\subset\overline{U}_n$ is the marked curve itself. So any
point $p\in\overline{U}_n$ belonging to the fibre over $C$
is actually a point on the stable $n$-pointed tree $C$, and may
therefore be considered as an
additional mark; stabilization may be necessary to ensure that
our new $(n+1)$-marked tree is stable. Vice-versa, given an $(n+1)$ pointed curve $C'$, we can think of the $(n+1)$st point as being a point on the $n$-marked curve obtained by forgetting the last marked point (eventually contracting, if needed); this way $C'$ corresponds to a point on the universal family $\overline{U_n}$. These constructions realize an isomorphism between
 $\overline{U}_n$ and $\overline{M}_{0,n+1}$. 

\subsubsection{The boundary}
We define the boundary to be the complement of $M_{0,n}$ in $\overline{M}_{0,n}$.
It consists of all nodal stable curves.  

\textbf{Fact:} the boundary is a union of irreducible components, 
corresponding to the different possible ways of arranging the marks on the various twigs; the codimension of a boundary component equals the number of nodes
in the curves in that component. See~\cite{k:kf} for more details.

The codimension $1$ boundary strata of
$\overline{M}_{0,n}$, called the \textbf{boundary divisors},
 are in one-to-one correspondence with all ways of
partitioning $[n]=A\cup B$ with the cardinality of both $ A$ and  $B$ strictly greater than 1.\\
A somewhat special class of boundary divisors  consists of those with only two 
marked points on a twig. Together, these components are sometimes called  the \textbf{soft boundary} and denoted
by $D_{i,j}$. We can think of $D_{i,j}$ as the image of the $i$th section, $\sigma_i$, of the $j$th forgetful map, $\pi_j$ (or vice-versa). 

\begin{figure}[htbp]
\begin{center}
\includegraphics[width=13cm]{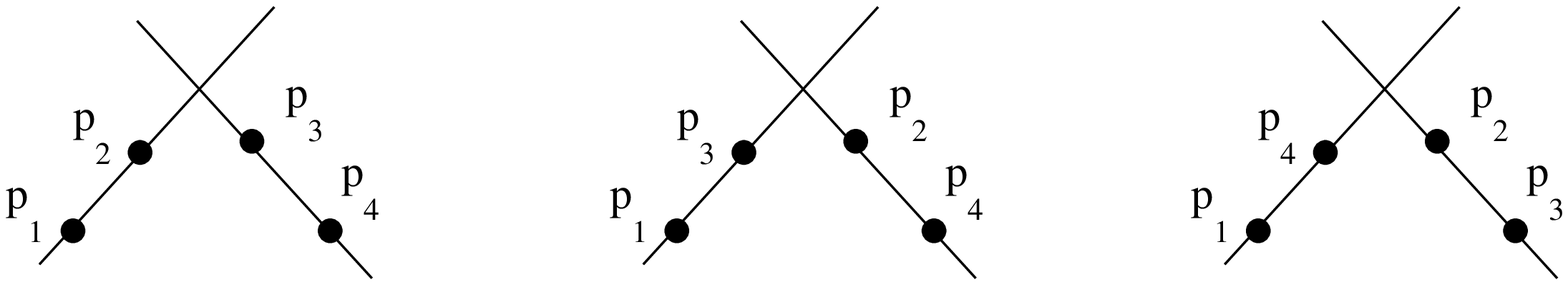}
\caption{irreducible components of the boundary of $\overline{M}_{0,4}$}\label{m04}
\end{center}
\end{figure}

\begin{figure}[htbp]
\begin{center}
\includegraphics[width=16cm]{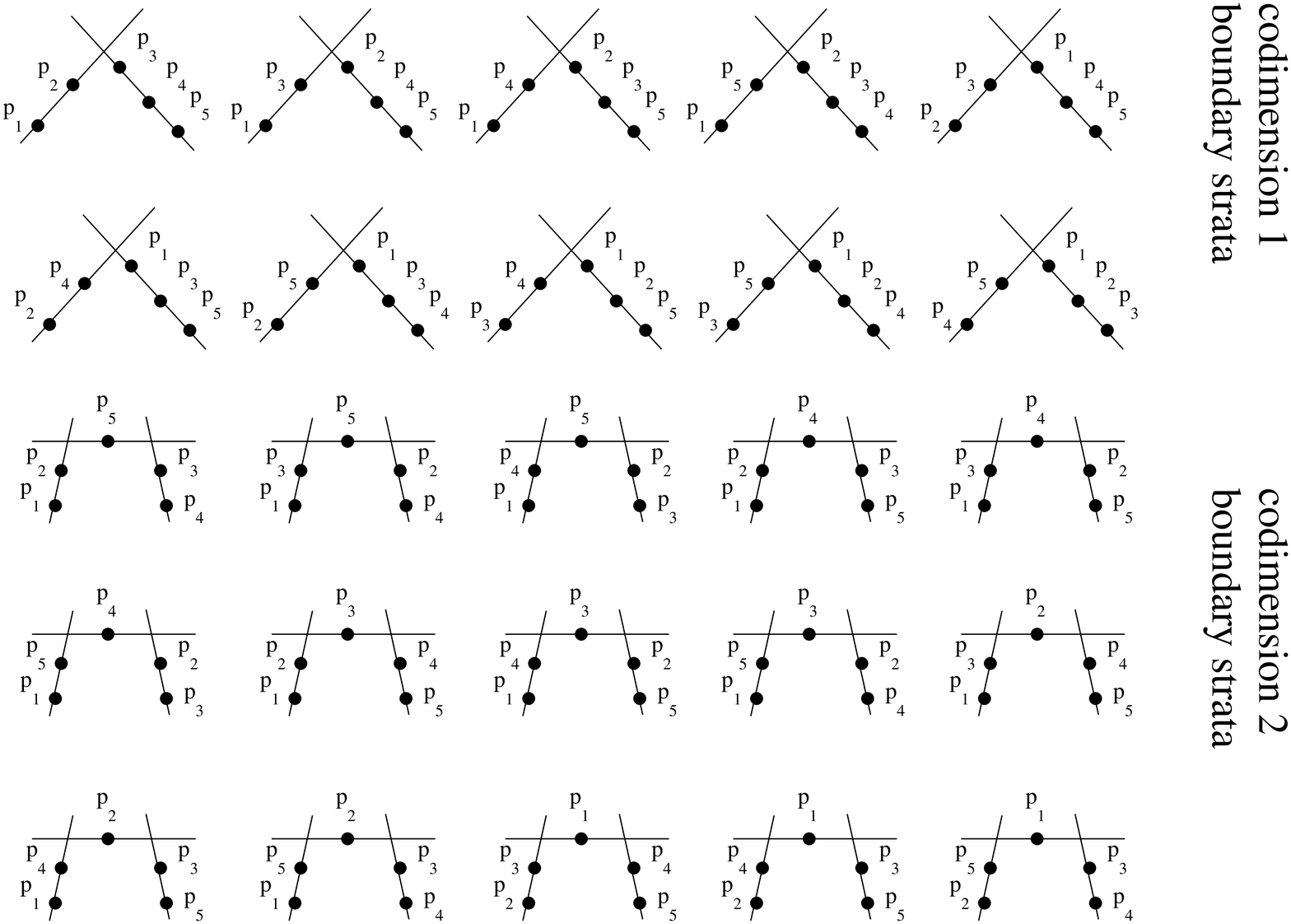}
\caption{boundary cycles of $\overline{M}_{0,5}$}\label{m05}
\end{center}
\end{figure}

There is
plenty more to be said about the spaces
$\overline{M}_{0,n}$, their relationships, and their boundaries, but
we will leave our treatment of $\overline{M}_{0,n}$ here, suggesting~\cite{k:kf}
as an excellent reference for beginners. In Figures~\ref{m04} and~\ref{m05} we draw all boundary strata for $\overline{M}_{0,4}$ and $\overline{M}_{0,5}$.

\subsection{Moduli of Rational Stable Maps}

Let us now move on to the moduli spaces of  greatest interest for
solving the bitangent problem. We would like to study, in general,
rational curves in projective space. The characteristic property of an
irreducible rational curve is that it can be parametrized by the
projective line, $\proj$. For this reason, it is natural to study maps
$\mu:\proj\rightarrow\mathbb{P}^r$. 

When we talk about the \textbf{degree} of such a map we mean the degree of
$\mu_\ast[\proj]$ as in homology. Be careful, the degree of the map may be different
from the degree of the image curve! For example the map
$$
\begin{array}{cccc}
\mu: &\proj      & \rightarrow & \Proj \\
     & (x_0:x_1) & \mapsto     & (x_0^2:x_1^2 :0) \\
\end{array}
$$
has degree two, but its image is a line.

Define $W(r,d)$ as the space of all
 maps from $\proj$ to $\mathbb{P}^r$
of degree $d$. A map in $W(r,d)$ is  specified, up to a constant, by $r+1$ binary
forms of degree
$d$ that do not all vanish at any point. It can then
be seen that $dim W(r,d)=(r+1)(d+1)-1$. 

We also have the family
$$
\begin{array}{ccc}
W(r,d)\times\proj & \stackrel{\rho}{\rightarrow} & \mathbb{P}^r\\
\downarrow &  &\\
W(r,d) & &
\end{array}
$$
where $\rho(\mu,x)=\mu(x)$. This family is tautological in the sense that the
fibre over the map $\mu$ is the map $\rho\mid_{\{\mu\}\times\proj}=\mu$. In fact, this
is a universal family. Thus, $W(r,d)$ is a fine moduli space for maps
$\proj\rightarrow\mathbb{P}^r$ of degree $d$. 

However, $W(r,d)$ is not the moduli space that we would like to study.
 For one, it is not compact. For another,
reparametrizations of the source curve are considered as different points in
$W(r,d)$. To fix the latter problem, let us simply consider the space
$M_{0,0}(\mathbb{P}^r,d):=W(r,d)/Aut(\proj)$.  For a detailed account on why this quotient is indeed a space, see chapter five of ~\cite{hm:mc}

 Another way to eliminate automorphisms is to consider $n$-pointed maps (maps
$\mu:C\rightarrow\mathbb{P}^r$  with an $n$-marked source
$C\simeq\proj)$. It should be no surprise that there is a fine
moduli space
$M_{0,n}(\mathbb{P}^r,d)$ for isomorphism classes of $n$-pointed maps
$\proj\rightarrow\mathbb{P}^r$ of degree $d$, namely $M_{0,n}\times
W(r,d)$. But we still have not dealt with the non-compactness of this
moduli space. The idea is to parallel the construction that led us to stable curves.
\begin{defi}
An \textbf{n-pointed stable map} is a map $\mu:C\rightarrow\mathbb{P}^r$, where:
	\begin{enumerate}
		\item $C$ is a $n$-marked tree. 
		\item Every twig in $C$ mapped to a point must have at least three
			 special points on it. 
	\end{enumerate}
\end{defi}

\textbf{Fact:} Moduli spaces of  $n$-pointed rational stable maps to $\mathbb{P}^r$ of degree $d$ (denoted $\overline{M}_{0,n}(\mathbb{P}^r,d))$ can be constructed; they compactify the moduli spaces of smooth maps.
It is
straightforward to verify that an $n$-pointed map is stable if and only if
it has only a finite number of automorphisms. Unfortunately, there is no way
to eliminate all nontrivial automorphism. Details can be found in ~\cite{k:kf}.

\textbf{Example:} An element
$\mu\in\overline{M}_{0,2}(\mathbb{P}^2,2)$  that is  the double cover of a
line, marking the ramification points,
admits a nontrivial automorphism  exchanging the two covers. 
This allows us to construct a nontrivial family of maps $\mu_t$ that maps constantly to one point in the moduli space. Consider:
$$
\begin{array}{cccc}
\mu_t: &[0,1]/\{0=1\} \times \proj & \longrightarrow & \Proj \\
 & & & \\
 & (t,(x_0:x_1)) & \mapsto & (0:x_0^2: e^{2\pi it}x_1^2).
\end{array}
$$  

Because of this phenomenon   there is no universal
family associated to the spaces
$\overline{M}_{0,n}(\mathbb{P}^r,d)$, and the corresponding moduli spaces are only coarse.

\vspace{0.5cm}
Since $M_{0,n}(\mathbb{P}^r,d)$ is dense in
$\overline{M}_{0,n}(\mathbb{P}^r,d)$, the latter has
dimension
$$(r+1)(d+1)-1+(n-3)=rd+r+d+n-3.
$$

\textbf{Example:} in particular, $\overline{M}_{0,n}(\mathbb{P}^r,0) = \overline{M}_{0,n} \times \mathbb{P}^r.$
\vspace{0.5cm}

There are some useful maps among the spaces
$\overline{M}_{0,n}(\mathbb{P}^r,d)$. For instance, as with the spaces
$\overline{M}_{0,n}$, we have the forgetful maps $\pi_i$
defined by simply forgetting the $i$th mark and the sections,
$\sigma_j$, of the family
$$
\begin{array}{c}
\overline{M}_{0,n+1}(\mathbb{P}^r,d)\\
\pi_i\downarrow\uparrow\sigma_j\\
\overline{M}_{0,n}(\mathbb{P}^r,d)
\end{array}
$$
defined by declaring the $j$th and  the $i$th mark to coincide. Contraction and stabilization 
are performed to make these maps defined everywhere. 

In addition,
there are 
\textbf{evaluation maps}
$$
\nu_i:\overline{M}_{0,n}(\mathbb{P}^r,d)\rightarrow\mathbb{P}^r
$$
defined by $\nu_i(\mu)=\mu(p_i)$ where $p_i$ is the $i$th mark on
the source curve $C$.

The forgetful and evaluation morphisms allow us to identify $\overline{M}_{0,n+1}(\mathbb{P}^r,d)$
as a tautological family for $\overline{M}_{0,n}(\mathbb{P}^r,d) $:

$$
\begin{array}{ccc}
\overline{M}_{0,n+1}(\mathbb{P}^r,d) & \stackrel{\nu_{n+1}}{\longrightarrow} &
\mathbb{P}^r\\
\pi_{n+1}\downarrow &  & \\
\overline{M}_{0,n}(\mathbb{P}^r,d). & &
\end{array}
$$

This way we can think of points of $\overline{M}_{0,n+1}(\mathbb{P}^r,d)$ either as $n+1$ pointed maps
or as points on an $n$-marked curve mapped to $\mathbb{P}^r$. Being comfortable with this identification
comes in very handy when making computations.

\subsubsection{The Boundary}
The boundary consists of maps
whose domains are reducible curves. In fact, its description  is very
similar to that of $\overline{M}_{0,n}$. Boundary strata are determined now not only
by the combinatorial data of the arrangement of the marks, but also by the degree  the
maps restrict to on each twig.

Boundary divisors are in
one to one correspondence with all ways of partitioning $[n]=A\cup B$
and $d=d_A+d_B$ such that:
\begin{itemize}
	\item $\#A\geq 2$ if $d_A=0$;
 	\item $\#B\geq 2$ if
$d_B=0$.
\end{itemize}

\subsection{Psi classes}

Consider a family of curves admitting a section.
$$
\begin{array}{c}
\mathcal{C}\\
\pi\downarrow\uparrow \sigma\\
B
\end{array}
$$
 We can define the \textbf{cotangent line bundle}, $\mathbb{L}_\sigma$, as the line bundle on $B$
whose fibre at a point $b\in B$ is the cotangent space of $\mathcal{C}_b=\pi^{-1}(b)$ at
 the point $\sigma(b)$. We call the $\mathbf{\psi}$ \textbf{class} of the family the first Chern class
of this line bundle.
Observe
 that, for a trivial family of curves 
with a constant section, the $\psi$ class vanishes.

This construction can be extended in a natural way to the moduli space 
$\overline{M}_{0,n}(\mathbb{P}^r,d)$.\\
Informally, we have a sheaf on the tautological family $\overline{M}_{0,n+1}(\mathbb{P}^r,d)$
 whose local sections  away from  nodes 
are differential forms
on the curves. We  obtain it by considering $1$-forms on the tautological family and quotienting 
by forms that are pulled back from the moduli space. This sheaf is called the \textbf{relative dualizing sheaf}\footnote{The word relative refers to the fact that we are quotienting by everything coming from downstairs. In other words, we are constructing a sheaf on the universal family of the moduli space by ``gluing'' together sheaves defined on the curves.} ,
and denoted by $\omega_{\pi_{n+1}}$.

\vspace{0.5cm}
Consider now the $i$th tautological section $\sigma_i$. If we restrict $\omega_{\pi_{n+1}}$ to this section,
we obtain a line bundle on the moduli space whose  fibres are naturally identified with the
cotangent spaces of the curves at the $i$th marked point. Then we can define the class:
$$
\mathbf{\psi}_i:=\mathbf{c}_1(\sigma_i^\ast\omega_{\pi_{n+1}})\in
A^1(\overline{M}_{0,n}(\mathbb{P}^r,d)).
$$

The construction of $\psi_i$ is natural in the sense that if we  have a family of pointed stable maps, inducing
a morphism to the moduli space, the   $\psi_i$ class of the family is the pull-back of the 
$\psi_i$ class on the moduli space.
\vspace{0.5cm}

It may seem that there is no difference between the information
carried by $\mathbf{\psi}_i\in
A^1(\overline{M}_{0,n+1}(\mathbb{P}^r,d))$ and $\mathbf{\psi}_i\in
A^1(\overline{M}_{0,n}(\mathbb{P}^r,d))$. We
have a natural map between these two spaces, the
tautological family
$\pi_{n+1}:\overline{M}_{0,n+1}(\mathbb{P}^r,d)\rightarrow\overline{M}_{0,n}(\mathbb{P}^r,d)$.
It may seem that
$\mathbb{L}_{i,n+1}:=\sigma_i^\ast\omega_{\pi_{n+1}}$ and
$\pi_{n+1}^\ast \mathbb{L}_{i,n}$ are the same line bundle, thus yielding the
relation $\mathbf{\psi}_i=\pi_{n+1}^\ast\mathbf{\psi}_i$.
In reality, this is almost true. Surely these line bundles agree off the component $D_{i,n+1}$
of the soft boundary. 
From this consideration, we conclude 
\begin{eqnarray}
\mathbb{L}_{i,n+1}=\pi_{n+1}^\ast \mathbb{L}_{i,n}\otimes\mathcal{O}(mD_{i,n+1})
\end{eqnarray}
for some integer $m$.

Next, observe that $\mathbb{L}_{i,n+1}$ restricted to $D_{i,n+1}$ is a
trivial line bundle: we are looking at curves with a twig having only
three special points; the node, the $i$th and the $(n+1)$st mark. By
an automorphism of the twig, we can assume that the node is at $0$ and
the two marks are at $1$ and $\infty$. Therefore, this line bundle
restricted to $D_{i,n+1}$ is the cotangent space at a single
unchanging point of $\proj$.  This implies
$$
 \mathcal{O}_{D_{i,n+1}} = {\mathbb{L}_{i,n+1}}_{\mid D_{i,n+1}} = 
$$
$$ = (\pi_{n+1}^\ast
\mathbb{L}_{i,n}\otimes\mathcal{O}(mD_{i,n+1}))_{\mid D_{i,n+1}} =
\mathbb{L}_{i,n} \otimes{\mathcal{O}(mD_{i,n+1})}_{\mid D_{i,n+1}}.$$

By the adjunction formula (\cite{gh:ag}, page 146), the line bundle \linebreak
$\mathcal{O}{(D_{i,n+1})}_{\mid D_{i,n+1}}$ is the normal bundle of the divisor $D_{i,n+1}$. 

But the normal 
directions to a section in the moduli space are precisely the  tangent directions to the fibres.
Hence  $\mathcal{O}{(D_{i,n+1})}_{\mid D_{i,n+1}}$
 is  the dual to the relative cotangent bundle $\mathbb{L}_{i,n}$. It follows that  $m$ must
be $1$.

Finally, by taking Chern classes in $(3)$, we can deduce the fundamental relation:
\begin{eqnarray}
\psi_i=\pi_{n+1}^\ast \psi_i + D_{i,n+1}.
\end{eqnarray}
The above pull-back relation can be used to describe explicitly $\psi$ classes for moduli spaces of rational stable curves
in terms of boundary strata. A closed formula can be found in~\cite{k:pc}.

\section{Counting Bitangents}
\subsection{The Strategy}
We now have all the necessary machinery to tackle our problem of counting bitangents.
Before we start digging deep into details and computations, let us outline our strategy.\\
	Let $Z:=\{f=0\}$ be a projective plane curve of degree $d$:
\begin{itemize}
	\item  we  consider the moduli 
space $\maps$, of two-pointed, rational stable maps of degree one;
	\item we   construct
a jet bundle on this space with the property that the zero set of a section of this bundle
consists of stable maps having at least second order contact with $Z$ at the image of the $i$th marked
point; we name this cycle $\Phi_i(Z)$;
	\item we  represent $\Phi_i(Z)$ in the Chow ring in terms of
$\psi$ classes and other natural classes;
	\item we step by step compute the intersection $\Phi_1(Z)\Phi_2(Z)$;
	\item we  identify and clean up some 
garbage that lives in that intersection and corresponds to maps that are not bitangents;
	\item finally, we have counted two pointed maps that are tangent to $Z$ at each mark; we just need
to divide by $2$, since we are not interested in the ordering of the marks.
\end{itemize}
Easy enough? Now let us start over slowly and do everything carefully.

\subsection{$\map$}

This space
has dimension 3, and it is explicitly realized by the following incidence relation:  
$$
\map
=
\{(p,l)
\in
\Proj\times\check{\Proj}
\mid
p\in
l\}
=:\I
\subset
\Proj\times\check{\Proj}.
$$
There are two projections of $\mathcal{I}$ onto $\Proj$ and $\check{\Proj}$,
 that we will denote $q$ and $\check{q}$.  The latter makes $\I$ into a
tautological family of lines in $\Proj$:
$$
\begin{array}{ccc}
\proj &
\longrightarrow &
\I \\
 & & \downarrow\check{q} \\
 & & \check{\Proj}.
\end{array}
$$
This family is tautological in the sense that the fibre over
$l\in\check{\Proj}$ is $l$ itself. 

A fibre over $p\in\Proj$ under
$q$ is the pencil of lines in $\Proj$
passing through $p$ and so this projection also gives rise to a
$\proj$ bundle over $\Proj$. Observe that $q$ is precisely the evaluation
map $\nu_x$. 

\textbf{Notation:} We denote by $x$ the unique mark in the space of one-pointed maps, and add the subscript $x$ to any entity (class, map ...) related to it. We do so to keep track of the conceptual difference from the marked points on the two-pointed maps, which will be numbered $1$ and $2$. 

\begin{defi}
We identify and name two natural divisors  on $\map$.
\begin{description}
	\item[$\mathbf{\iota}(p):$] Let us look at 
the  hyperplane divisor $\check{p}\subset\check{\Proj}$
of all lines
passing through a point $p$, and consider the cycle of its pull-back
$\check{q}^\ast(\check{p}):=\mathbf{\iota}(p)$.
	\item[$\mathbf{\eta}_x(l):$] Similarly, we look at the hyperplane divisor
$l\subset\Proj$  and define $\mathbf{\eta}_x(l):=q^\ast (l)$ as  its pull-back
under $q$. 
\end{description}
\end{defi}

In general, define
$\mathbf{\eta}_x(Z):=q^\ast(Z)=\nu^\ast_x(Z)$ as the cycle of maps whose
mark is sent into
$Z$. 

There is only one class of a line and
only one of a point in $A^\ast(\Proj)$, hence
$\mathbf{\eta}_x:=[\mathbf{\eta}_x(l)]$ and
$\mathbf{\iota}:=[\mathbf{\iota}(p)]$ are independent of $l$ and $p$
respectively.\\
 Since $\map$ is a $\proj$ bundle over $\Proj$, it follows that
$\mathbf{\iota}$ and $\mathbf{\eta}_x$, i.e. the pull-backs of hyperplane divisors in the base and in the fiber, generate $A^1(\map)$. It
is therefore useful to know all intersections of the two
classes.

It is a good exercise to construct a picture verifying each
of the following statements about classes.
\begin{itemize}
	\item$\mathbf{\eta}_x=[\{(p^\prime,l^\prime)\mid p^\prime\in l,
 l^\prime\in\check{p^\prime}$ with $l$ fixed$\}]$
	\item$\mathbf{\iota}=[\{(p^\prime,l^\prime)\mid l^\prime\in\check{p}, 
  p^\prime\in l^\prime$ with $p$ fixed$\}]$
\item
$\mathbf{\eta}_x^2=[\{(p,l^\prime)\mid l^\prime\in\check{p}$ with $p$
  fixed$\}]$
\item
$\mathbf{\iota\eta}_x=[\{(p^\prime,l^\prime)\mid p^\prime\in l,
  l\in\check{p}$ with $p$ and $l$ fixed$\}]$
\item
$\mathbf{\iota}^2=[\{(p^\prime,l)\mid p^\prime\in\l$ with $l$ fixed$\}]$
\item
$\mathbf{\iota\eta}_x^2=[\{(p,l^\prime)\mid l\in\check{p}\cap\check{p}_0$
  with $p,p_0$ fixed$\}]=\mathbf{pt}$ since there is exactly one line
  passing through two distinct points.
\item
$\mathbf{\eta}_x\iota^2=[\{(p^\prime,l)\mid p^\prime\in l\cap l_0$ with
  $l,l_0$ fixed$\}]=\mathbf{pt}$ since there is exactly one point in
  the intersections of two distinct lines.
\item
$\mathbf{\eta}_x^3=[\{(p^\prime,l^\prime)\mid
  l^\prime\in\check{p}_1\cap\check{p}_2\cap\check{p}_3, p^\prime\in
  l^\prime$ with $p_1,p_2,p_3$ fixed$\}]=0$ since in general three
  points do not lie on a common line.
\item
$\mathbf{\iota}^3=[\{(p^\prime,l^\prime)\mid p^\prime\in l_1\cap
    l_2\cap l_3, p^\prime\in l^\prime$ with $l_1,l_2,l_3$ fixed$\}]=0$
    since in general three
    lines do not share a common point.
\end{itemize}

The following two lemmas prove identities that will be crucial for our later computations.
\begin{lemma}
\begin{eqnarray}
\mathbf{\eta}_x\iota=\mathbf{\eta}_x^2+\mathbf{\iota}^2. \label{etaiota}
\end{eqnarray}
\end{lemma}
\textsc{Proof:}
We  construct a
one-parameter family of cycles, parametrized by $[0,1]$, with the
left-hand side of our identity
as one endpoint of this family
and the right-hand side as the other.\\ 
To choose a representative of the class 
$\mathbf{\eta}_x\iota$, one must specify a fixed point
$p$ and a fixed line $l$. Let us fix $l$ once and for
all and let $p_t$ be a path in $\Proj$ such that
$p_t\in l$ if and only if $t=0$.  Our one parameter family
$\mathbf{\alpha}_t$ is defined as follows: 
$$\mathbf{\alpha}_t=\{(p^\prime,l^\prime)\mid p^\prime\in l, l\in\check{p_t}\}.$$
 Notice that
$[\mathbf{\alpha}_t]=\mathbf{\eta}_x\iota$ for $t\not=0$.
 
\begin{figure}[htbp]
\begin{center}
\includegraphics[width=12cm]{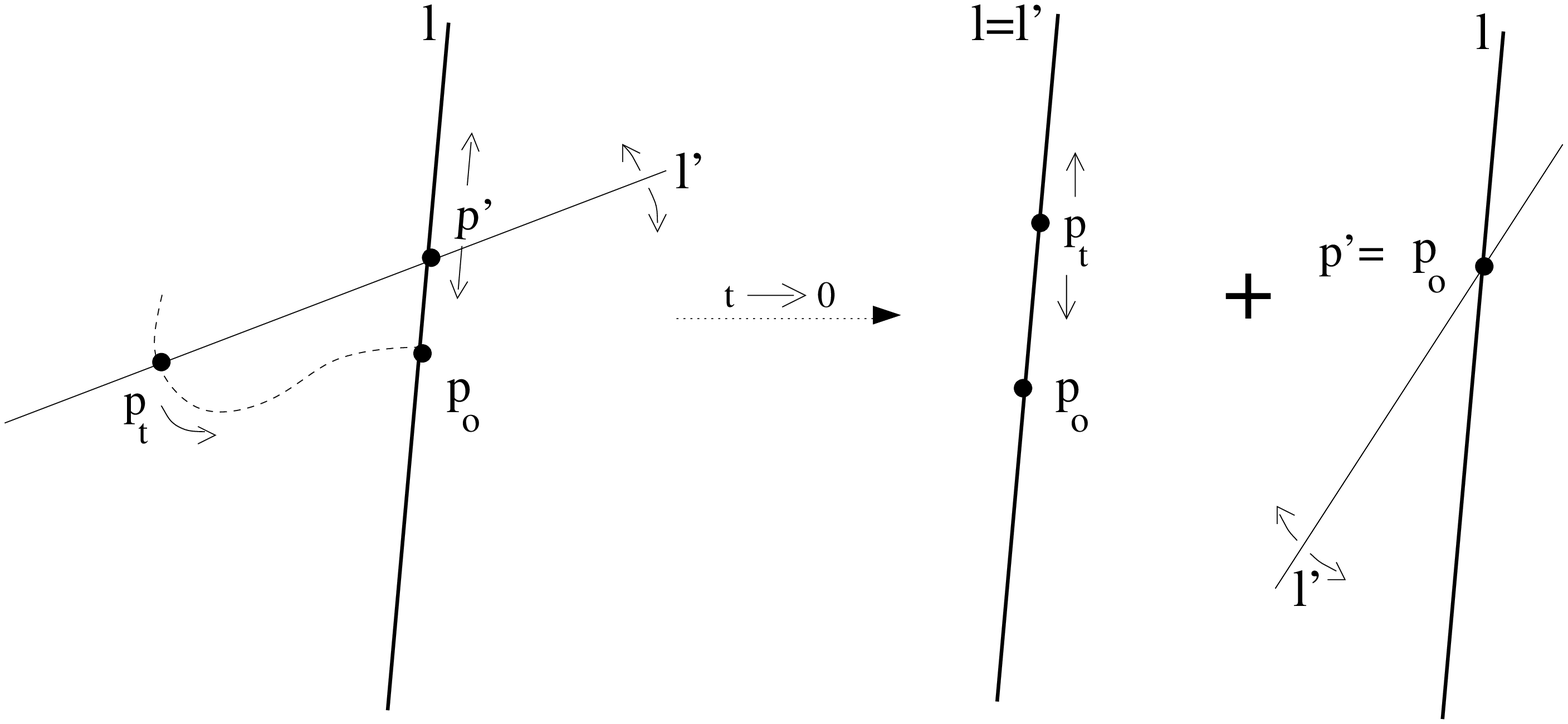}
\caption{$\mathbf{\eta}_x\iota=\mathbf{\eta}_x^2+\mathbf{\iota}^2.$}
 \end{center}
\end{figure}

 We  now
examine what happens as $t\rightarrow0$.
\begin{description}
\item[$p^\prime\not=p_0$:]
 necessarily, $l^\prime=l$ and our resulting one dimensional class
is parametrized solely by $p^\prime\in l$ and is therefore
$\mathbf{\iota}^2$.
\item [$p^\prime=p_0$:]
Then $l^\prime$ is only required to be in $\check{p_0}$ and so such
$l^\prime$'s in $\check{p_0}$ parametrize our resulting one
dimensional class. We have arrived at $\mathbf{\eta}_x^2$.
\end{description}

We now have that $\alpha_0=\eta_x^2+\iota^2$, which allows us to conclude (\ref{etaiota}).
\vspace{0.5cm}

\begin{lemma}
\begin{eqnarray}
\mathbf{\psi}_x=\mathbf{\iota}-2\mathbf{\eta}_x. \label{psi}
\end{eqnarray}
\end{lemma}
\textsc{Proof:} As $\mathbf{\psi}_x\in A^1(\map)$, it is possible to express
$\mathbf{\psi}_x=a\mathbf{\iota}+b\mathbf{\eta}_x$ for some integers $a$
and $b$. Let us determine these two integers.
\begin{description}
	\item[a:] intersecting $\psi_x$ with $\eta_x^2$ we obtain
$$\psi_x\eta_x^2=a\mathbf{\iota}\mathbf{\eta}_x^2+b\mathbf{\eta}_x^3=a.$$
 
Consider
$\sigma_x^\ast\omega_{\pi_2}$ restricted to $\eta_x^2=\{(p,l^\prime)\mid p$
is fixed and $l^\prime\in\check{p}\}$. It is the line
bundle over $\eta_x^2$ whose fibre over a point $(p,l^\prime)\in\eta_x^2$ is the
cotangent space of $l^\prime$ at the fixed point $p$.

It is worth convincing yourself that this is
the dual of the tautological line bundle
$$
\begin{array}{cl}
\mathcal{S} & \\
\downarrow\pi & \\
\proj & = \{l\subset\Proj\mid p\in l\} = \eta_x^2.
\end{array}
$$

We computed (section \ref{Chern}) that $\mathbf{c}_1(\mathcal{S})=-1$
and so $a=1$.

	\item[b:]
similarly,
$b$ is the product $\mathbf{\psi}_x{\iota}^2$. To find this
intersection we must consider the line bundle
$\sigma_x^\ast\omega_{\pi}$ restricted to $\iota^2=\{p\in l\mid l$ is
fixed$\}$. A fibre of this line bundle over a point $p\in l$ is the
cotangent space of our fixed $l$ at $p$. This is simply the
cotangent bundle of $l$.

Since $l=\proj=S^2$ has Euler
characteristic $2$, then the degree of the first Chern class of the cotangent bundle 
is $-2$. We thus have that $b=-2$.
\end{description}

\subsection{$\maps$}
 First off note that $dim\maps=4$. The
description of $\maps$ is slightly more
complicated largely due to the existence of its only boundary divisor which we
 call $\beta$.

\begin{figure}[htbp]
\begin{center}
\includegraphics[height=4cm]{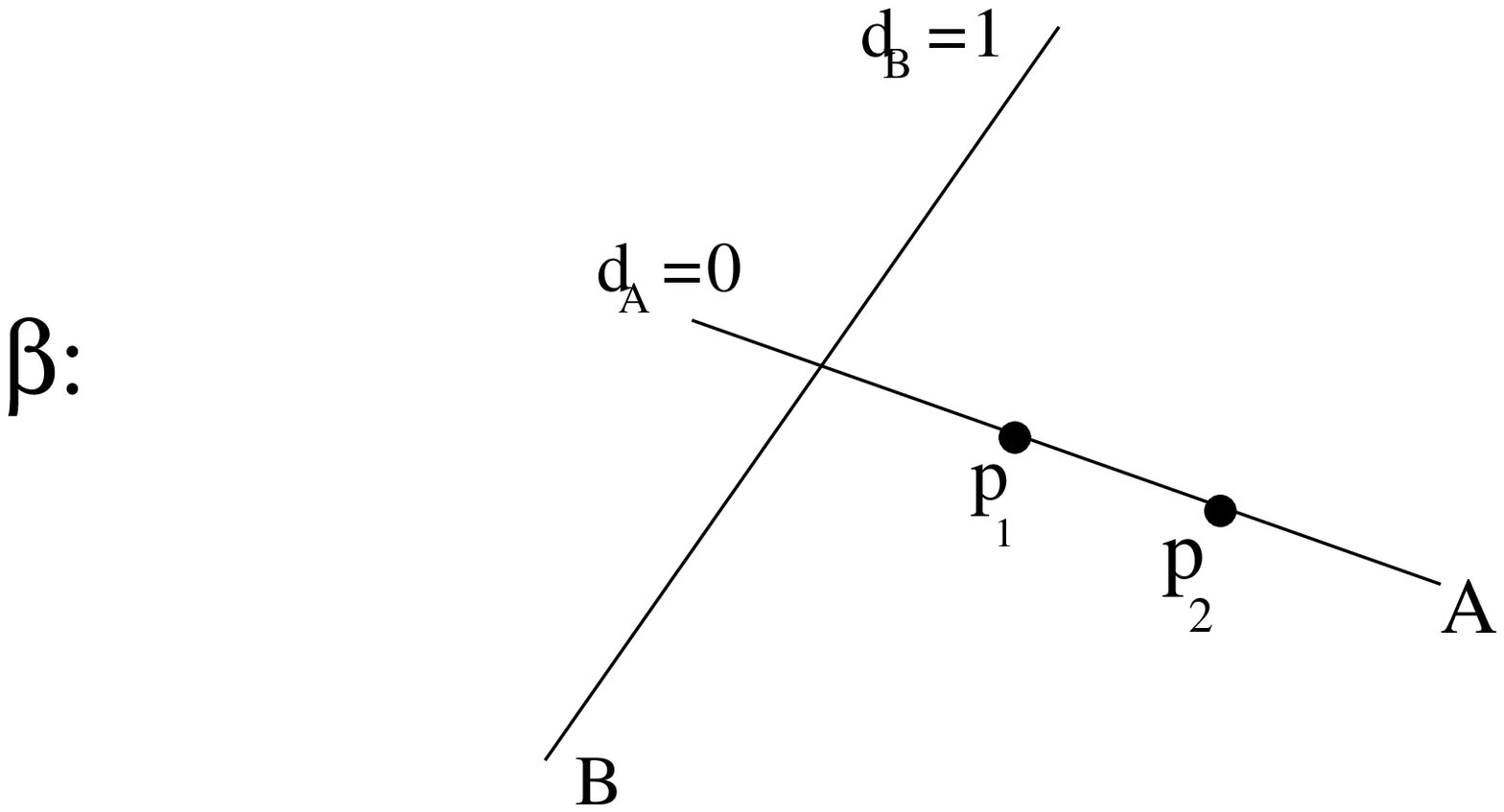}
\caption{the boundary divisor in $\maps$}
 \end{center}
\end{figure}

For a  two-pointed map $\mu$ not in the boundary, $\mu(p_1)\not=\mu(p_2)$. Since we
are considering maps of degree 1, i.e. isomorphisms of lines, $\mu(p_1)$ and
$\mu(p_2)$ completely determine (up to reparametrization)j our map $\mu$. It follows that
$$\maps\setminus\beta=\Proj\times\Proj\setminus\Delta.$$
 On $\beta$, $\mu(p_1)=\mu(p_2)$: for any
line $l$ through $p$, there is a map in $\beta$ that contracts the twig with the marks to
$p$ and maps the other twig isomorphically to $l$.\\
So for our
description to be complete, we need to replace
$(p,p)\in\Delta\subset\Proj\times\Proj$ with a $\proj$ worth of
maps. We arrive at
$$
\maps=Bl_\Delta(\Proj\times\Proj).
$$

Consider the tautological families
$$
\begin{array}{c}
\maps\\
\pi_1\downarrow\downarrow\pi_2\\
\map.
\end{array}
$$
Both families have a natural common section $\sigma:=\sigma_1=\sigma_2$. The image of 
$\sigma$ is the unique boundary divisor $\beta$ in
$\maps$.

Define $\eta_i(Z):= \nu_i^\ast(Z)$ It's easy to check the following identities:
\begin{enumerate}
	\item$\sigma_\ast(\mathbf{\eta}_x(Z))=\mathbf{\beta\eta}_1(Z)=\mathbf{\beta\eta}_2(Z)$.
	  \item $\sigma^\ast(\eta_i(Z))=\eta$.
	\item$\pi_i^\ast\mathbf{\eta}_x(Z)=\mathbf{\eta}_i(Z)$.
\end{enumerate}

\subsection{Tangency Conditions}
Let us define
$\mathbf{\Phi}_i(Z)\in
A^\ast(\overline{M}_{0,n}(\mathbb{P}^2,1))$ as
the cycle of maps  tangent   to a plane curve $Z=\{f=0\}$ at the image of the $i$th marked point.  Formally,
$$\mathbf{\Phi}_i(Z)=\{\mu\in\overline{M}_{0,n}(\mathbb{P}^2,1)|\mu^\ast
f \mbox{ vanishes at }p_i\mbox{ with
multiplicity }\geq2\}.$$ 

We now want to obtain an expression for $\mathbf{\Phi}_x(Z)\in
A^\ast(\map)$ in terms of $\mathbf{\eta}_x$, $\mathbf{\beta}$, and
$\mathbf{\psi}_x$.\\
Consider the tautological family
$$
\begin{array}{ccc}
\maps & \stackrel{\nu_2}{\longrightarrow}  & \Proj    \\
 & & \\
\sigma\uparrow \downarrow \pi_2& \nearrow \nu_x &    \\
 &  & \\
\map. & & 
\end{array}
$$
Let us pull-back the line bundle $\mathcal{O}(Z)$ via $\nu_2$, and consider 
the  first jet bundle $J_{\pi_2}^1\nu_2^\ast\mathcal{O}(Z)$ relative to $\pi_2$. Relative here
means that we quotient out by everything that can be pulled back from $\map$. 
Let $Z$ be defined
by the vanishing of the polynomial $f$, and let us consider                 the
 zero locus of the section
$\tau:=\nu_2^\ast f + (\partial_{\pi_2}^1\nu_2^\ast f) dt\in \Gamma(J_{\pi_2}^1\nu_2^\ast\mathcal{O}(Z)$); what we 
obtain is the  
 locus of maps in $\maps$ such that the pull-back of $f$ at the second marked point vanishes to second order. If we interpret $\maps$ as the universal family for $\map$ \footnote{this is true because one pointed maps of degree one have no nontrivial automorphisms!}, it follows that to obtain $\mathbf{\Phi}_x(Z)$, the locus in $\map$ of
lines tangent to $Z$ at the unique marked point, we need to pull-back via the section $\sigma$.
 
In formulas, this translates to
 $$\mathbf{\Phi}_x(Z) =\mathbf{e}(\sigma^\ast
J_{\pi_2}^1\nu_2^\ast\mathcal{O}(Z)).$$
Since the rank of the bundle in question is $2$, the Euler class is the second Chern class.
\vspace{0.5cm}

 To calculate $\mathbf{c}_2(\sigma^\ast
J_{\pi_2}^1\nu_2^\ast\mathcal{O}(Z))$ we  use the
following short exact sequence discussed in section \ref{jets}.
$$
0\rightarrow\nu_2^\ast\mathcal{O}(Z)\otimes\omega_{\pi_2}\rightarrow
J_{\pi_2}^1\nu_2^\ast\mathcal{O}(Z)\rightarrow\nu_2^\ast\mathcal{O}(Z)\rightarrow0.$$
Notice that the first and last terms of this sequence are line
bundles. We then want to consider the pull-back along $\sigma$ of
this exact sequence. Using the Whitney formula, we now have that
$$
\mathbf{\Phi}_x(Z)=\mathbf{c}_1(\sigma^\ast\nu_2^\ast\mathcal{O}(Z))\mathbf{c}_1(\sigma^\ast\nu_2^\ast\mathcal{O}(Z)\otimes\sigma^\ast\omega_{\pi_2})=d\mathbf{\eta}_x(d\mathbf{\eta}_x+\mathbf{\psi}_x)
$$ 
in $\map$.

The last equality follows from the two facts:
\begin{itemize}
  \item $[Z]= d\mathbf{H} \in A^1(\Proj)$, where $\mathbf{H}$ is the hyperplane class generating $A^\ast(\Proj)$ and $d=deg f$;
    \item $\nu_x = \nu_2\sigma$, and $\eta_x$ is by definition $\nu_x^\ast(\mathbf{H})$.
\end{itemize}
\vspace{0.5cm}

Now we want to consider the case when we have more than one marked point: let us say we want to compute $\Phi_1(Z)$ in $\maps$. The obvious guess is 
$\Phi_1(Z)=d\mathbf{\eta}_1(d\mathbf{\eta}_1+\mathbf{\psi}_1)$.
We need to be careful, though:
for
maps in $\beta\mathbf{\eta}_1(Z)\subset d\mathbf{\eta}_1(d\mathbf{\eta}_1+\mathbf{\psi}_1)$, the whole twig containing the
two marks is mapped to $Z$. So, for $\mu$ such a map, $\mu^\ast f$
vanishes identically along the contracting twig and thus to all
orders. We do not want to consider these maps as tangents to $Z$.

\textbf{Fact:}
This simple correction works. The formula in $\maps$ is
\begin{eqnarray}
\mathbf{\Phi}_1(Z)=d\mathbf{\eta}_1(d\mathbf{\eta}_1+\mathbf{\psi}_1-\mathbf{\beta}).
\end{eqnarray}
Lastly, note that
$\mathbf{\beta\Phi}_1(Z)=\sigma_\ast\mathbf{\Phi}_x(Z)$.

\vspace{0.5cm} 
\textbf{Remark:} It would be nice to use higher order jet bundles to
describe cycles of maps having higher order contact with our curve
$Z$. Unfortunately in general it is quite difficult, as fairly
complicated problems of excess intersection arise. 

For our
application, we only need to push our luck a little further: we need
to describe the cycle $\mathbf{\Phi}_x^{(3)}(Z)$ of inflection
tangents to $Z$.  Luckily, thanks to the fact that $\map$ has no
boundary, the argument carries through: $\mathbf{\Phi}_x^{(3)}(Z)$ can
be computed as the Euler class of the bundle $E^\prime:=\sigma^\ast
J_{\pi_2}^2\nu_2^\ast\mathcal{O}(Z)$. Here, rank$ E^\prime=3$. By the exact sequence (\ref{jet}) and the Whitney formula:
$$\mathbf{\Phi}_x^{(3)}(Z)=\mathbf{c}_3(E^\prime)=\mathbf{\eta}_x(Z)(\mathbf{\eta}_x(Z)+\mathbf{\psi}_x)(\mathbf{\eta}_x(Z)+2\mathbf{\psi}_x).
$$
Using our previous calculations,
\begin{eqnarray}
\mathbf{\Phi}_x^{(3)}(Z)=(3d^2-6d)\mathbf{pt}.
\end{eqnarray}
We will abuse notation from now on and leave off the class of a point
in our calculations. When writing a dimension zero class we will
simply write its integral over the fundamental class.

\subsection{The Computation}
Let us finally get down to business. We define in
general
 $$\Lambda_Z(m_1p_1+m_2p_2)\subset\maps$$
 as the cycle of maps
$\mu$, such that $\mu^\ast\mathbf{Z}\geq m_1p_1+m_2p_2$. Note that $\mathbf{\Lambda}_Z(2p_1)=\mathbf{\Phi}_1(Z)$.\\
 
Our ultimate goal is to compute
$\mathbf{\Lambda}_Z(2p_1+2p_2)$, i.e. the class of 
  maps in $\maps$ which are tangent to $Z$ at both $p_1$ and $p_2$.\\
 Note that
  $dim\mathbf{\Lambda}_Z(2p_1+2p_2)=0$ since we are
  imposing 4 independent conditions in
  a space of dimension 4. This tells us that our enumerative problem makes sense.

On first thought, one might suggest
$\mathbf{\Lambda}_Z(2p_1+2p_2)=\mathbf{\Phi}_1(Z)\mathbf{\Phi}_2(Z)$.
After all, $\mathbf{\Phi}_1(Z)$ is the class of all maps tangent to
$Z$ at $p_1$ and similarly for $\mathbf{\Phi}_2(Z)$ at $p_2$ so their
intersection
seems to be what we are after. However, one must
be  careful: for example, $\mathbf{\beta\Phi}_1(Z)$ is in this
intersection and we don't want to consider such maps.

We will proceed in two steps and obtain $\mathbf{\Lambda}_Z(2p_1+2p_2)$ from
$\mathbf{\Phi}_1(Z)\mathbf{\Phi}_2(Z)$ by "throwing away the spurious maps". 
\vspace{0.5cm}

\textbf{Step 1:}
We consider  the intersection
$\mathbf{\Phi}_1(Z)\mathbf{\eta}_2(Z)$, consisting  of all maps tangent to $Z$ at
$p_1$ that intersect
$Z$ at $p_2$. What we get are two parts: 
\begin{itemize}
  \item $\mathbf{\Lambda}_Z(2p_1+p_2)$, in which $p_1$ and $p_2$
do not lie on the same twig of degree zero;
\item
$\mathbf{\beta\Phi}_1(Z)=\sigma_\ast\mathbf{\Phi}_x(Z)$, corresponding to maps with both marked points on a degree $0$ twig.
\end{itemize}
Set theoretically,
$$
\eta_2(Z)\Phi_1(Z)=\Lambda_Z(2p_1+p_2)\cup\sigma_\ast\Phi_x(Z).
$$
The correct multiplicities are 1 and 2
respectively, giving:
$$
\mathbf{\Lambda}_Z(2p_1+p_2)=\mathbf{\eta}_2(Z)\mathbf{\Phi}_1(Z)-2\sigma_\ast\mathbf{\Phi}_x(Z).
$$

\textbf{Step 2:}We now intersect $
\mathbf{\Lambda}_Z(2p_1+p_2)$ with
$\mathbf{\eta}_2(Z)+\mathbf{\psi}_2-\mathbf{\beta}$. This imposes second order vanishing at the second marked point.

Again, this intersection gives one part with multiplicity 1, which is
not on the boundary $\mathbf{\beta}$, and another part with
multiplicity 2 in $\mathbf{\beta}$. The first is what we are looking
for: $\mathbf{\Lambda}_Z(2p_1+2p_2)$. To find the other part, remember
we are already
working inside $\mathbf{\Lambda}_Z(2p_1+p_2)$, so to lie in
  $\mathbf{\beta}$ means $p_2\rightarrow p_1$. We thus obtain
  $\mathbf{\beta\Lambda}_Z(3p_1)=\mathbf{\beta\Phi}_1^{(3)}$ with
  multiplicity 2.

\textbf{Note:} To explain these multiplicities rigorously is a subtle business
  which we will not go into here. The intuition behind these $2$'s is that the boundary contribution can be carried by either of the classes we are intersecting, and hence shows up twice in the intersection. The reader interested in how to carry out these computation can consult \cite{g:tnopc}.

Finally, we have identified

\begin{eqnarray}
\mathbf{\Lambda}_Z(2p_1+2p_2) &=& -2\sigma_\ast\mathbf{\Phi}_x^{(3)}(Z)+(z\mathbf{\eta}_2+\mathbf{\psi}_2-\mathbf{\beta})\mathbf{\Lambda}_Z(2p_1+p_2)\nonumber\\
&=& -2\sigma_\ast\mathbf{\Phi}_x^{(3)}+\mathbf{\Phi}_1(Z)\mathbf{\Phi}_2(Z)-2\sigma_\ast\mathbf{\Phi}_x(Z)(z\mathbf{\eta}_2+\mathbf{\psi}_2-\mathbf{\beta}).\nonumber \\
\end{eqnarray}

Now all that stands in our way of calculating bitangents are
substitutions and computations.

\vspace{0.5cm}

Before we move on, take a second to remember or derive the following easy facts, that we will use in the forthcoming computations:
\begin{description}\label{list}
  \item [(a)]$\mathbf{\iota\eta}_x^2=\mathbf{\iota}^2\mathbf{\eta}_x=1$
    \item [(b)]$\mathbf{\iota\eta}_x=\mathbf{\eta}_x^2+\mathbf{\iota}^2$
      \item [(c)]$\mathbf{\psi}_x=\mathbf{\iota}-2\mathbf{\eta}_x$ 
        \item [(d)]$\mathbf{\eta}_i^3=0$
	  \item [(e)]$\mathbf{\psi}_1=\pi_2^\ast\mathbf{\psi}_x+\mathbf{\beta}$
	    \item [(f)]$\mathbf{\beta\psi}_i=0$
	      \item [(g)]$\mathbf{\beta}^2=-\mathbf{\beta}\pi_i^\ast\mathbf{\psi}_x$
	        \item [(h)]$\mathbf{\beta\eta}_1=\mathbf{\beta\eta}_2$
		  \item [(i)]$\sigma_\ast\mathbf{\Phi}_x(Z)=\mathbf{\beta\Phi}_1(Z)=\mathbf{\beta\Phi}_2(Z)$
		     \item [(j)]$\mathbf{\eta}_1\pi_2^\ast\mathbf{\alpha}=\pi_2^\ast(\mathbf{\alpha\eta}_x)$
for any class $\mathbf{\alpha}\in A^\ast(\map)$.
\end{description}

Let us expand the first term in (9),  $\sigma_\ast\mathbf{\Phi}_x^{(3)}$.
We have already found $\mathbf{\Phi}_x^{(3)}(Z)=(3d^2-6d)$ in $\map$. Since
pushing forward preserves dimension and $\sigma_i$ is
injective, thus sending the class of a
point to the class of a point, then our first term becomes 
$$\sigma_\ast\mathbf{\Phi}_x^{(3)}=3d^2-6d.$$ 

Now for the second and third terms.
\begin{eqnarray}
\mathbf{\Phi}_1(Z)\mathbf{\Phi}_2(Z) & =
& d\mathbf{\eta}_1(d\mathbf{\eta}_1+\mathbf{\psi}_1-\mathbf{\beta})d\mathbf{\eta}_2(d\mathbf{\eta}_2+\mathbf{\psi}_2-\mathbf{\beta})\nonumber\\
 & = &
d^4(\mathbf{\eta}_1^2\mathbf{\eta}_2^2)+d^3(\mathbf{\eta}_1^2\mathbf{\eta}_2\mathbf{\psi}_2+\mathbf{\eta}_1\mathbf{\eta}_2^2\mathbf{\psi}_1)+d^2(\mathbf{\eta}_1\mathbf{\eta}_2\mathbf{\psi}_1\mathbf{\psi}_2+\mathbf{\eta}_1^2\mathbf{\beta}^2)\nonumber.
\end{eqnarray}
\begin{eqnarray}
\sigma_\ast\mathbf{\Phi}_x(Z)(d\mathbf{\eta}_2+\mathbf{\psi}_2-\mathbf{\beta})
& = &
\mathbf{\beta}((d\mathbf{\eta}_x)(d\mathbf{\eta}_1+\mathbf{\psi}_1-\mathbf{\beta}))(d\mathbf{\eta}_2+\mathbf{\psi}_2-\mathbf{\beta})\nonumber\\
& = & -2d^2(\mathbf{\eta}_1^2\mathbf{\beta}^2)+d(\mathbf{\eta}_1\mathbf{\beta}^3)\nonumber.
\end{eqnarray}
We must now compute each intersection in the above expressions\footnote{The little numbers over the equal signs refer to the identities from page \pageref{list} that are used at each step.}. 
\begin{itemize}
  \item$\mathbf{\eta}_1^2\mathbf{\eta}_2^2=1$, as seen from the figure \ref{eta} illustrating the fact that there is exactly one line passing
through two prescribed points.

\begin{figure}
\begin{center}
\includegraphics[height=5cm]{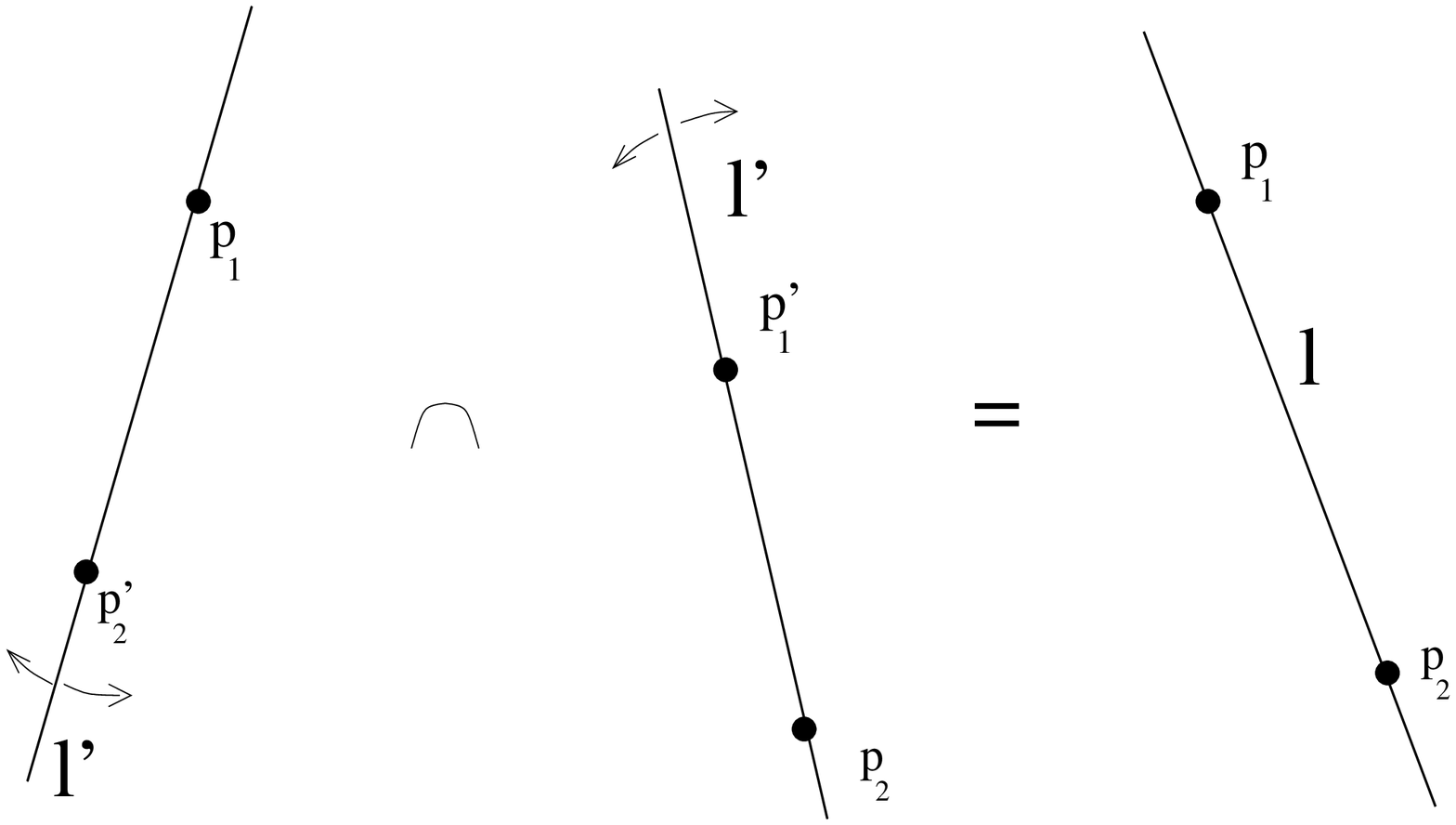}
\caption{the intersection $\mathbf{\eta}_1^2\mathbf{\eta}_2^2$. }
\label{eta}
 \end{center}
\end{figure}
\item$\mathbf{\eta}_1^2\mathbf{\beta}^2
\ \stackrel{(f)}{=}\ \mathbf{\eta}_1^2\mathbf{\beta}(-\pi_2^\ast\mathbf{\psi}_x)
\ \stackrel{(i)}{=}\ -\mathbf{\beta}\pi_2^\ast(\mathbf{\eta}_x^2\mathbf{\psi}_x)
\ \stackrel{(b+c)}{=}\ -\mathbf{\beta}\pi_2^\ast(a)=-\mathbf{\beta}[fibre]$.\\
But $\beta$ is the image of a section, hence it
intersects all fibres transversely. Thus,
$\mathbf{\eta}_1^2\mathbf{\beta}^2=-1$.
\item$\mathbf{\eta}_1\mathbf{\beta}^3
\ \stackrel{(f)}{=}\ -\mathbf{\eta}_1\mathbf{\beta}^2\pi_2^\ast\mathbf{\psi}_x
\ \stackrel{(f)}{=}\ \mathbf{\eta}_1\mathbf{\beta}(\pi_2^\ast\mathbf{\psi}_x)^2
=\mathbf{\eta}_1\mathbf{\beta}(\pi_2^\ast\mathbf{\psi}_x^2)
\ \stackrel{(i)}{=}\ \mathbf{\beta}\pi_2^\ast(\mathbf{\psi}_x^2\mathbf{\eta}_x)
\ \stackrel{(a+b+c)}{=}\ \mathbf{\beta}\pi_2^\ast(-3)=-3$.
\item$\mathbf{\eta}_1^2\mathbf{\eta}_2\mathbf{\psi}_2
\ \stackrel{(d)}{=}\mathbf{\eta}_1^2\mathbf{\eta}_2(\pi_1^\ast\mathbf{\psi}_x+\mathbf{\beta})
\ \stackrel{(i)}{=}\mathbf{\eta}_1^2\pi_1^\ast(\mathbf{\eta}_x\mathbf{\psi}_x)
\ \stackrel{(a+b)}{=}\mathbf{\eta}_1^2\pi_1^\ast(\mathbf{\iota}^2-\mathbf{\eta}_x^2)$.\\
 Notice first of all that $\mathbf{\eta}_1^2\pi_1^\ast\mathbf{\eta}_x^2= \mathbf{\eta}_1^2\mathbf{\eta}_2^2=1$, as shown in figure \ref{eta}. 

Next, we claim that
$\mathbf{\eta}_1^2\pi_1^\ast\mathbf{\iota}^2=0$. In fact $\mathbf{\eta}_1^2$ is the class of all two-pointed lines passing
through a fixed point, where the first mark is at the fixed point while the second mark is free to move.
\begin{figure}[htbp]
\begin{center}
\includegraphics[height=5cm]{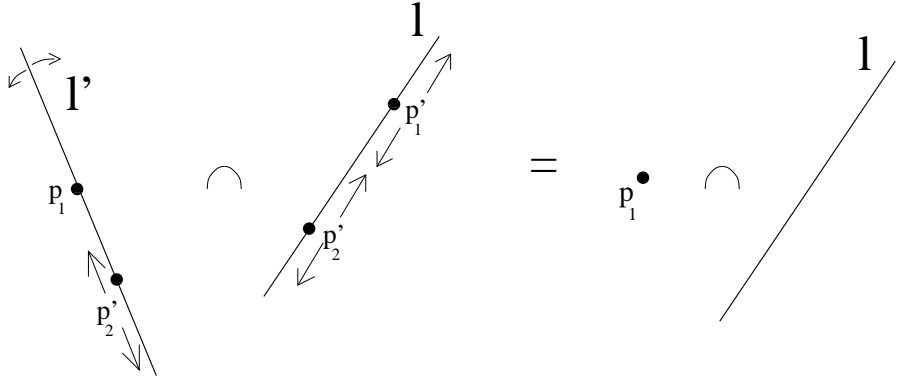}
\caption{the intersection $\mathbf{\eta}_1^2\pi_1^\ast(\mathbf{\iota}^2)$. }
\label{neta}
 \end{center}
\end{figure}
 
Intersecting  with $\pi_1^\ast(\mathbf{\iota}^2)$ means to require that
our fixed point intersects our fixed line transversely, which is to say
that they don't intersect at all. This is illustrated in figure \ref{neta}.

 Thus,
$\mathbf{\eta}_1^2\mathbf{\eta}_2\mathbf{\psi}_2=-1$; by 
symmetry, we also have
$\mathbf{\eta}_2^2\mathbf{\eta}_1\mathbf{\psi}_1=-1$.

\item Finally,
\begin{eqnarray}
\hspace{-1.5cm}\mathbf{\eta}_1\mathbf{\eta}_2\mathbf{\psi}_1\mathbf{\psi}_2
& \stackrel{(d)}{=} &
\mathbf{\eta}_1(\pi_2^\ast\mathbf{\psi}_x+\mathbf{\beta})\mathbf{\eta}_2(\pi_1^\ast\mathbf{\psi}_x+\mathbf{\beta})\nonumber\\
& \stackrel{(i)}{=} &
(\pi_2^\ast(\mathbf{\psi}_x\mathbf{\eta}_x)+\mathbf{\beta\eta}_1)(\pi_1^\ast(\mathbf{\psi}_x\mathbf{\eta}_x+\mathbf{\beta\eta}_2)\nonumber\\
& \stackrel{(a+b)}{=} &
\pi_1^\ast(\mathbf{\iota}^2-\mathbf{\eta}_x^2)\pi_2^\ast(\mathbf{\iota}^2-\mathbf{\eta}_x^2)+\mathbf{\beta\eta}_2\pi_1^\ast(\mathbf{\iota}^2-\mathbf{\eta}_x^2)+\mathbf{\beta\eta}_1\pi_2^\ast(\mathbf{\iota}^2-\mathbf{\eta}_x^2)+\mathbf{\beta}^2\mathbf{\eta}_1^2\nonumber\\
& = &
\pi_1^\ast(\mathbf{\iota}^2)\pi_2^\ast(\mathbf{\iota}^2)-\pi_1^\ast(\mathbf{\iota}^2)\pi_2^\ast(\mathbf{\eta}_x^2)-\pi_1^\ast(\mathbf{\eta}_x^2)\pi_2^\ast(\mathbf{\iota}^2)+\pi_1^\ast(\mathbf{\eta}_x^2)\pi_2^\ast(\mathbf{\eta}_x^2)\nonumber\\
& &+\mathbf{\beta}\pi_1^\ast(\mathbf{\iota}^2\mathbf{\eta}_x)+\mathbf{\beta}\pi_2^\ast(\mathbf{\iota}^2\mathbf{\eta}_x)+\mathbf{\eta}_x^2\mathbf{\beta}^2.\nonumber
\end{eqnarray}
\end{itemize}
Via our previous calculations, the last three terms can be easily seen
as $1,1,$ and $-1$ respectively. So we must now find the first four
terms above. We  do so by recalling pictures.

\vspace{0.5cm}
$\pi_1^\ast(\mathbf{\iota}^2)$ is the
class of a fixed line and all ordered pairs of points on it.
 
To
intersect two such classes is to require that our line is fixed as 
\textit{two}
transverse lines, which is impossible. Thus
$\pi_1^\ast(\mathbf{\iota}^2)\pi_2^\ast(\mathbf{\iota}^2)=0$.

\vspace{0.5cm}
$\pi_1^\ast(\mathbf{\eta}_x^2)= \eta_2^2$, and notice that now symmetry implies that $\eta_2^2\pi_2^\ast(\mathbf{\iota}^2)=\eta_1^2\pi_1^\ast(\mathbf{\iota}^2)= 0$, as shown in figure \ref{neta}. Similarly,
$\pi_2^\ast(\mathbf{\eta}_x^2)\pi_1^\ast(\mathbf{\iota}^2)=0$.

\vspace{0.5cm}
Lastly, let us 
intersect $\pi_2^\ast(\mathbf{\eta}_x^2)$ with
$\pi_1^\ast(\mathbf{\eta}_x^2)$. But it's clear that
$\pi_2^\ast(\mathbf{\eta}_x^2)\pi_1^\ast(\mathbf{\eta}_x^2)=\eta_1^2\eta_2^2$, and we have shown in figure \ref{eta} that this intersection is $1$.

We have then found
$$
\mathbf{\eta}_1\mathbf{\eta}_2\mathbf{\psi}_1\mathbf{\psi}_2=0+0+0+1+1+1-1=2.
$$

Putting this all together, we have now calculated
\begin{eqnarray}
\mathbf{\Lambda}_Z(2p_1+2p_2) & = &
-2(3d^2-6d)+(d^4-2d^3+d^2)-2(-2d^2-3d)\nonumber\\
& = & d^4-2d^3-9d^2+18d\nonumber\\
& = & d(d-2)(d-3)(d+3).\nonumber
\end{eqnarray}

After the dust has settled, we now know that
 a generic plane curve  $Z$  of degree $d$ has
$$
N_{\mathcal{B}}(d)=\frac{1}{2}d^4-d^3-\frac{9}{2}d^2+9d
$$
bitangents. Remember we are dividing by $2$ because we do not care about the order of the marked points.

Notice that for $d=2$ and $d=3$ we get that there are no bitangents as should 
be the case. For $d=4$ we find $28$ bitangents, 
the first interesting result.

%
%
%
\bibliographystyle{alpha}
\bibliography{renzo}

%






%
%
%
\Addresses
%
%
%
%
\end{document}